\documentstyle[12pt]{article}
\newcommand{\np} {} 

\newcommand{\parf}{\subsection}
\newcommand{\punk}{\subsubsection}

\newcommand{\bold}{\bf}
\font\bfit  = cmbxti10 scaled \magstep1
 
\newcommand{\rbox} [1]{{\rm{#1}}}        
\newcommand{\bbox} [1]{{\bf #1}}          
\newcommand{\ttbox}[1]{\rbox{\tt #1}}       
\expandafter\chardef\csname pre amsym.def 
at\endcsname=\the\catcode`\@
\catcode`\@=11
\def\newsymbol#1#2#3#4#5{\let\next@\relax
\ifnum#2=\@ne\let\next@\msafam@\else
 \ifnum#2=\tw@\let\next@\msbfam@\fi\fi
 \mathchardef#1="#3\next@#4#5}
\def\hexnumber@#1{\ifcase#1 0\or 1\or 2\or 3\or 4\or 5\or 6\or 
7\or 8\or 9\or A\or B\or C\or D\or E\or F\fi}
\newfam\msafam
\newfam\msbfam
\edef\msafam@{\hexnumber@\msafam}
\edef\msbfam@{\hexnumber@\msbfam}
\def\Bbb#1{\fam\msbfam\relax#1}
\catcode`\@=\csname pre amsym.def at\endcsname
\newsymbol\subsetneqq  2324
\newsymbol\restriction 1316
\newsymbol\complement  107B
\newsymbol\preceQ      1334
\newsymbol\succeQ      133C
\newsymbol\boX         1003
\font\twlmsa=msam10 scaled \magstep1
\font\tenmsa=msam10 
\font\egtmsa=msam8
\def\xiimsa{\textfont\msafam=\twlmsa\scriptfont\msafam=\egtmsa}
\def\xmsa  {\textfont\msafam=\tenmsa\scriptfont\msafam=\egtmsa}
\font\twlmsb=msbm10 scaled \magstep1
\font\tenmsb=msbm10
\font\egtmsb=msbm8
\def\xiimsb{\textfont\msbfam=\twlmsb\scriptfont\msbfam=\egtmsb}
\def\xmsb  {\textfont\msbfam=\tenmsb\scriptfont\msbfam=\egtmsb}
\font\twleufm=eufm10 scaled \magstep1
\font\teneufm=eufm10
\font\egteufm=eufm8
\newfam\eufmfam
\newcommand{\xiieufm}{\textfont\eufmfam=\twleufm
    \scriptfont\eufmfam=\egteufm\def\fraK{\fam\eufmfam\twleufm}}
\newcommand{\xeufm}{\textfont\eufmfam=\teneufm
\def\fraK{\fam\eufmfam\teneufm}}
\newcommand{\xiieusm}{\font\teusm=eusm10 scaled \magstep1
                     \font\seusm=eusm8 \font\zeusm=eusm6}
\newcommand{\xeusm}{\font\teusm=eusm10 
                     \font\seusm=eusm7 \font\zeusm=eusm5}
\newcommand{\skr}[1]{{\mathchoice {\hbox{\teusm{#1}}} 
{\hbox{\teusm{#1}}}{\hbox{\seusm{#1}}}{\hbox{\zeusm{#1}}} }}
\def\addto#1#2{
\ifx\zzone\undefined\let\zzone=#1\def#1{\zzone#2}\else
\ifx\zztwo\undefined\let\zztwo=#1\def#1{\zztwo#2}\else
\fi\fi
}
\addto\normalsize   {\xiimsb\xiimsa\xiieufm\xiieusm} 
\addto\footnotesize {\xmsb\xmsa\xeufm\xeusm} 
\newtheorem{theorem}             {Theorem}
\newtheorem{corollary}  [theorem]{Corollary}
\newtheorem{definition} [theorem]{Definition}
\newtheorem{lemma}      [theorem]{Lemma}
\newtheorem{proposition}[theorem]{Proposition}
\newtheorem{remark}     [theorem]{Remark}
\newtheorem{asser} {{\bfit Assertion}}         
\newtheorem{aq}    {{\bfit Acknowledgements}}  
\newtheorem{bla}   {{\bfit Blanket assumption}}
\newtheorem{prF}   {Proof}                     
\newtheorem{prFi}  {{\bfit Proof}}             
\newcommand{\thsp}{\hspace{0.3ex}}
\newcommand{\bass}{\begin{asser}\thsp} 
\newcommand{\eass}{\end{asser}} 
\newcommand{\bcor}{\begin{corollary}\thsp}
\newcommand{\ecor}{\end{corollary}}
\newcommand{\bdf} {\begin{definition}\rm\thsp}
\newcommand{\edf} {\end{definition}}
\newcommand{\ble} {\begin{lemma}\thsp}
\newcommand{\ele} {\end{lemma}}
\newcommand{\bre} {\begin{remark}\rm\thsp}
\newcommand{\ere} {\end{remark}}
\newcommand{\bte} {\begin{theorem}\thsp}
\newcommand{\ete} {\end{theorem}}
\newcommand{\bpro}{\begin{proposition}\thsp} 
\newcommand{\epro}{\end{proposition}} 
\newcommand{\bbla}{\begin{bla}\rm\thsp}
\newcommand{\ebla}{\end{bla}}
\newcommand{\baq} {\begin{aq}\rm\thsp}
\newcommand{\eaq} {\end{aq}}
\newcommand{\bpf} {\begin{prF}\rm} 
\newcommand{\epf} {\qed\end{prF}} 
\newcommand{\bpfi}{\begin{prFi}\rm} 
\newcommand{\epfi}{\qedi\end{prFi}} 
\newcommand{\epF} [1]{\qeD{#1}\end{prF}} 
\newcommand{\qed}   {\hfill$\mtho\boX$} 
\newcommand{\qedi}  {\hfill$\mtho\dashv$} 
\newcommand{\qeD}[1]{\hfill\hbox{$\mtho\boX$ 
({\it#1\hspace{0.3ex}\/})}} 
\newcommand{\ben}{\begin{enumerate}
        \itemsep=1mm plus 0.5mm minus 0.5mm}
\newcommand{\een}{\end{enumerate}}
\newcommand{\bit}{\begin{itemize}
        \itemsep=1mm plus 0.5mm minus 0.5mm}
\newcommand{\eit}{\end{itemize}}
\newcommand{\bay}{\begin{array}}
\newcommand{\eay}{\end{array}}
\newcommand{\ZFC} {\bbox{ZFC}}
\newcommand{\bP}  {\bbox{P}}
\newcommand{\rF}  {\rbox{F}}
\newcommand{\rH}  {\rbox{H}}
\newcommand{\rha} {{\rH^\ast}}
\newcommand{\rL}  {\rbox{L}}
\newcommand{\rU}  {\rbox{U}}
\newcommand{\rV}  {\rbox{V}}
\newcommand{\wo}  {\rbox{\tt WO}}
\newcommand{\hc}  {\rbox{HC}}

\newcommand{\Ord} {\ttbox{Ord}}
\newcommand{\dosp}{\hspace{0.5ex}}
\newcommand{\disp}{\hspace{0.4ex}}
\newcommand{\ea}  {\hbox{\it e.\disp a.\/}}
\newcommand{\ie}  {\hbox{\it i.\disp e.\/}}
\newcommand{\etc} {{\it etc.\/}}
\newcommand{\eg}  {\hbox{\it e.\disp g.\/}}
\newcommand{\po}  {\hbox{p.\dosp o.}}
\newcommand{\lo}  {\hbox{l.\dosp o.}}
\newcommand{\pqo} {\hbox{p.\dosp q.-o.}}

\newcommand{\hop} {\hbox{h.\dosp o.\dosp p.}}

\newcommand{\wlg} {\hbox{w.\dosp l.\dosp o.\dosp g.}}
\newcommand{\kL} {\rL}
\newcommand{\kV} {\rV}
\newcommand{\kvp}{{\kV^+}}

\newcommand{\al} {\alpha} 
\newcommand{\ba} {\beta} 
\newcommand{\ga} {\gamma} 

\newcommand{\Da} {\Delta}
\newcommand{\La} {\Lambda} 
\newcommand{\la} {\lambda} 
\newcommand{\kpa}{\kappa}
   \newcommand{\kpp}{{\kpa^+}}
\newcommand{\vpi}{\varphi}
\newcommand{\vvi}{{\vec\vpi}}
\newcommand{\vt} {\vartheta}
\newcommand{\sg} {\sigma}
\newcommand{\Sg} {\Sigma}
\newcommand{\om} {\omega} 
\newcommand{\omi}{\om_1}

\newcommand{\lom} {^{<\om}} 
\newcommand{\lomi}{^{<\omi}} 
\newcommand{\Ups}{\Upsilon}
\newcommand{\ups}{\upsilon}

\newcommand{\Rho}{\rbox{R}}
\newcommand{\Tau}{{\cal T}} 
\newcommand{\rdi} {\la}
\newcommand{\rdii}{(\rdi\pone)}

\newcommand{\pone}{\hspace{-0.4ex}+\hspace{-0.4ex}1}

\newcommand{\fs}[2]{{\bold\Sigma}^{#1}_{#2}}
\newcommand{\fp}[2]{{\bold\Pi}^{#1}_{#2}}
\newcommand{\fd}[2]{{\bold\Delta}^{#1}_{#2}}
\newcommand{\iSg}{{\mathchar"7106}}
\newcommand{\iPi}{{\mathchar"7105}}
\newcommand{\iDa}{{\mathchar"7101}}
\newcommand{\is}[2]{\iSg^{#1}_{#2}}
\newcommand{\ip}[2]{\iPi^{#1}_{#2}}
\newcommand{\id}[2]{\iDa^{#1}_{#2}}
\newcommand{\bbg}{\hspace{0.1ex}} 
\newcommand{\pr} {{\bbg{\fraK p}\bbg}\hspace{1pt}}
\newcommand{\gG} {{\fraK G}}
\newcommand{\gh} {{\bbox h}}
\newcommand{\gP} {{\fraK P}}
\newcommand{\gp} {{\fraK p}}
\newcommand{\bbd}{\hspace{0.1ex}}
\newcommand{\dvo}[1]{{\bbd{\Bbb#1}\bbd}}

\newcommand{\dP}    {{\bP}}
\newcommand{\dQ}    {\dvo Q}
\newcommand{\dZ}    {\dvo Z}
\newcommand{\dpw}   {\dP_{(2)}}
\newcommand{\dpwt}  {{\mathord{{}^T\hspace{-0.1ex}\dpw}}}
\newcommand{\dpws}  {{\mathord{{}^S\hspace{-0.1ex}\dpw}}}
\newcommand{\dpwo}  {\dP\tit\dP}

\newcommand{\skrsp}{\hspace{0.2ex}}
\newcommand{\skri}[1]{{\skrsp\skr{#1}\skrsp}}
\newcommand{\cF} {{\skri F}}
\newcommand{\cH} {{\skri H}}
\newcommand{\cL} {{\cal L}}
\newcommand{\cN} {{\skri N}}

\newcommand{\lan} [1]{\cL_{#1,0}}

\newcommand{\lano}   {\lan{\rdi+1}}

\newcommand{\kaf} {\cF}

\newcommand{\kah} {\cH}

\newcommand{\emps}{\emptyset}
\newcommand{\sq}  {\subseteq}

\newcommand{\cj}  {\mathbin{\hspace{0.2ex}\wedge\hspace{0.2ex}}}
\newcommand{\orr} {\mathbin{\textstyle\bigvee}}
\newcommand{\eqv} {\mathbin{\,\Longleftrightarrow\,}}
\newcommand{\imp} {\mathbin{\,\Longrightarrow\,}}
\newcommand{\lra} {\longrightarrow} 
\newcommand{\we}  {{\mathbin{\hspace*{0.2ex}^\wedge}}}
\newcommand{\sus} {{\exists\,}}
\newcommand{\kaz} {{\forall\,}}
\newcommand{\res} {{\hspace{0.2ex}\restriction\hspace{0.2ex}}}

\newcommand{\<}   {\leq}
\newcommand{\ti}  {\times}

\newcommand{\tit} {\mathbin{\ti_{T}}}
\newcommand{\dm}  {$$}
\newcommand{\mni} {\cei} 

\newcommand{\ima}{\mathbin{\hbox{\rm''}}}
\newcommand{\ang} [1] {\langle #1\rangle}
\newcommand{\cei} [1] {\lceil #1\rceil}
\newcommand{\ans} [1] {\{\hspace{0.1ex}#1\hspace{0.1ex}\}}
\newcommand{\stk}[2] {{\ang{#1\hspace{1pt};\hspace{1pt}#2}}}
\newcommand{\dop} [1] {\mathopen\complement{{#1}}}

\newcommand{\bve} {{\textstyle\bigvee\hspace{0.1ex}}}
\newcommand{\bwe} {{\textstyle\bigwedge\hspace{0.1ex}}}

\newcommand{\eqdf}{=_{\rbox{df}}}
\newcommand{\net} {\neg\,}
\newcommand{\mek}  {\preceQ}
   
\newcommand{\mekt} {\mek_T}
\newcommand{\meks} {\preceQ_S}
\newcommand{\nmek} {\not\mek}   
\newcommand{\nmekt}{\not\mekt}
\newcommand{\nmeks}{\not\meks}
\newcommand{\mekta} [1] {\mathrel{{\mek}_T^{#1}}}

\newcommand{\mektf} {\mek_{T,\dof}}

\newcommand{\meksf} {\mek_{S,\dof}}
\newcommand{\neksf} {\not\mek_{S,\dof}}

\newcommand{\mektF} {\mek_{T,f}}
\newcommand{\meksF} {\mek_{S,f}}

\newcommand{\Rmekt}[1] {\mathrel{{\mek}_{T,\,{#1}}}}
\newcommand{\Rmeks}[1] {\mathrel{{\mek}_{S,\,{#1}}}}
\newcommand{\meh} {\leq_\rH}
\newcommand{\eqh} {\equiv_\rH}

\newcommand{\eee} {\approx}
\newcommand{\eet} {\eee_T} 
\newcommand{\eesf}{\eee_{S,\dof}}
\newcommand{\nesf}{\not\eee_{S,\dof}}
\newcommand{\ees} {\eee_S}

\newcommand{\mso} {\mathbin{<_0}} 
\newcommand{\meo} {\mathbin{\leq_0}}
 
\newcommand{\mel} {\mathbin{\leq_{\rbox{lex}}}}
\newcommand{\hull}[3]
{{#1}_{#2}\hspace{0.25ex}[\hspace{0.5pt}#3\hspace{0.5pt}]}
\newcommand{\hua}[2]{\hull{\kaf}{#1}{#2}}

\newcommand{\huh}[2]{\hull{\kah}{#1}{#2}}
\newcommand{\hax}[1]  {\hua\dox{#1}}
\newcommand{\hay}[1]  {\hua\doy{#1}}
\newcommand{\haz}[1]  {\hua\doz{#1}}
\newcommand{\haxp}[1] {\hua\doxp{#1}}

\newcommand{\hazp}[1] {\hua\dozp{#1}}
\newcommand{\hhx}[1]  {\huh\dox{#1}}
\newcommand{\hhy}[1]  {\huh\doy{#1}}
\newcommand{\hhz}[1]  {\huh\doz{#1}}
\newcommand{\hhxp}[1] {\huh\doxp{#1}}

\newcommand{\itla}{\item\label}
\newcommand{\dd}[1]{$\mtho\hspace{0.2ex}{#1}$-}
\newcommand{\dl}  {\dd\la}

\font\ess=cmss8
\newcommand{\refo}[1]{{\mathchoice {\hbox{\sf{#1}}} 
{\hbox{\sf{#1}}}{\hbox{\ess{#1}}}{\hbox{\ess{#1}}} }}
\newcommand{\relf}[1]{{\refo #1}}
\newcommand{\qE}  {\mathbin{\relf E}}
\newcommand{\Eo}  {\mathbin{{\qE}_0}}
\newcommand{\nEo} {\mathbin{\not{{\hspace{-0.4ex}\qE}}_0}}

\newcommand{\dof} {{\dot f}}
\newcommand{\dtg} {{\dot g}}

\newcommand{\dox} {{\dot x}}
\newcommand{\doy} {{\dot y}}
\newcommand{\doz} {{\dot z}}
\newcommand{\doxp}{{\dot x'}}

\newcommand{\dozp}{{\dot z'}}
\newcounter{enuF}
\newcommand{\mtho}{\mathsurround=0mm}
\newcommand{\msur}{\hspace*{-1\mathsurround}}

\newcommand{\noi}{\noindent}
\newcommand{\vom}{\vspace{1mm}}
\newcommand{\vtm}{\vspace{2mm}}

\begin{document}

\normalsize

\title{Linearization of analytic order relations} 

\author{Vladimir Kanovei
\thanks{\ Moscow Transport Engineering Institute}
\thanks{\ {\tt kanovei@mech.math.msu.su} \ and \ 
{\tt kanovei@math.uni-wuppertal.de}
}
\thanks{\ This paper was accomplished in part during my visit to 
Caltech in April 1997. I thank Caltech for the support and 
A. S. Kechris and J. Zapletal for useful information and 
interesting discussions relevant to the topic of this paper 
during the visit.}
}
\date{April 1997} 
\maketitle
\normalsize

\begin{abstract}
We prove that if $\mek$ is an analytic partial order then either 
$\mek$ can be extended to a $\fd12$ linear order similar to an 
antichain in $2\lomi$ ordered lexicographically or a certain Borel 
partial order $\meo$ embeds in $\mek.$ Some corollaries for 
analytic equivalence relations are given, for instance, if $\qE$ 
is a $\is11[z]$ equivalence relation such that $\Eo$ does not 
embed in $\qE$ then $\qE$ is determined by intersections with 
\dd\qE invariand Borel sets coded in $\kL[z]$. 
\vtm \vtm\vtm

\hspace{0mm}\phantom{1.} Introduction\dotfill\pageref{int}\vom

\ref{prel}. Preliminaries\dotfill\pageref{prel}\vom

\ref{di}. Monotone Borel functions and the dichotomy\dotfill
\pageref{di}\vom

\ref f. The basic forcing\dotfill\pageref{f}\vom

\ref{pf}. The product forcing\dotfill\pageref{pf}\vom

\ref{E-b}. 
The construction of an embedding\dotfill\pageref{E-b}\vom

\ref{abs}. Why embedding $\meo$ is absolute\dotfill
\pageref{abs}\vom

\ref A. 
Borel and analytic order relations\dotfill\pageref{A}\vom

\ref S. Special cases: Borel classes and generic models\dotfill
\pageref{S}\vom

\hspace{0mm}\phantom{1.} References.\dotfill\pageref{refe}

\end{abstract}

\newpage

\parf*{Introduction}
\label{int}

%
It is a simple application of Zorn's lemma that any partial order 
can be extended to a linear order on the same domain. More 
generally any partial quasi-order admits a linearization.
\footnote
{\ {\bf Notation.\/} 
Several notions related to orders are sometimes understood 
differently, so let us take a space to fix an unambiguous meaning. 

A binary relation $\mek$ on a set $X$ is a 
{\it partial quasi-order\/}, or {\it\pqo\/} in brief, on $X,$ iff 
$x\mek y\cj y\mek z\imp x\mek z,$ and $x\mek x$ for any $x\in X.$ 
In this case, $\eee$ is the associated equivalence relation, 
\ie\ $x\eee y$ iff $x\mek y\cj y\mek x$. 

If $x\eee x\imp x=x$ for any $x$ then $\mek$ is 
a {\it partial order\/}, or {\it\po\/}. 
If in addition $x\mek y\orr y\mek x$ for all $x,\,y\in X$ then 
$\mek$ is a {\it linear\/} order (\lo).

Let $\mek$ and $\mek'$ be \pqo's on resp.\ $X$ and $X'.$ A map 
$h:X \lra X'$ will be called {\it half order preserving\/}, 
or {\it\hop\/}, iff 
$x\mek y\imp h(x)\mek' h(y)$.

Finally a {\it linearization\/} is any \hop\ map 
$h:\stk X\mek\lra\stk{X'}{\mek'},$ where $\mek'$ is a \lo, 
satisfying $x\eee y\eqv h(x)=h(y)$.
\vtm
\hfill $\mtho\bullet$
}

A much more difficult problem is to provide a descriptive 
characterization of the linear order in the assumption that one 
has such for the given \pqo. For instance, not every Borel \pqo\  
is {\it Borel linearizable\/}. \vtm

{\bfit Example 1.\/} Recall that $\Eo$ is an equivalence relation 
on $2^\om$ defined as follows: $a\Eo b$ iff $a(k)=b(k)$ for all 
but finite $k.$ Let $\mek$ be $\Eo$ considered as a \pqo. Then 
$\mek$ is not Borel linearizable. Indeed any linearization $h$ 
satisfies $a\Eo b\eqv h(a)=h(b),$ but it is known that Borel maps 
$h$ with such a property do not exist (see 
Harrington, Kechris, Louveau~\cite{hkl}).\vtm\hfill$\mtho\dashv$

{\bfit Example 2.\/} 
Example 1 can be converted to a partial order. Define the 
{\it anti-lexicographical\/} \po\ $\meo$ on $2^\om$ as follows:
$a\meo b$ iff either $a=b$ or there is $m\in\om$ such that 
$a(k)=b(k)$ for all $k>m$ and $a(m)<b(m).$ Clearly $a\meo b$ 
implies $a\Eo b$ and $\meo$ linearly orders each 
\dd{\Eo}equivalence class similarly to the integers $\dZ,$ except 
for the class of $\om\ti\ans{0}$ (ordered as $\om$) and 
the class of $\om\ti\ans{1}$ (ordered as $\om^{-1}$ -- the 
inverted $\om$)~\footnote
{\ If one enlarges $\mso$ so that, in addition, 
$a\mso b$ whenever $a,\,b\in 2^\om$ are such that $a(k)=1$ and 
$b(k)=0$ for all but finite $k$ then the enlarged relation can 
be induced by a Borel action of $\dZ$ on $2^\om,$ such that 
$a\mso b$ iff $a=zb$ for some $z\in \dZ,\msur$ $z>0$.} 
Finally $\meo$ is not Borel linearizable 
(see Subsection~\ref{bor}).\vtm\hfill$\mtho\dashv$

There are Borel-non-linearizable Borel orders of different nature, 
\eg\ the \pqo\ $a\mek b$ iff $a(k)\leq b(k)$ for all but finite 
$k$ on $2^\om$ or the dominance relation on $\om^\om.$ However by 
the next theorem the relation $\meo$ of Example 2 is actually a 
{\it minimal\/} Borel-non-linearizable Borel order. 
(Compare with the ``Glimm--Effros'' theorem of Harrington, Kechris, 
Louveau~\cite{hkl} saying that $\Eo$ is a minimal non-smooth 
Borel equivalence relation.)

\bte
\label{hb}
{\rm(Kanovei~\cite{k-b})} \ 
Suppose that\/ ${\mek}$ is a Borel \pqo\/ on\/ $\cN=\om^\om.$ 
Then exactly one of the following two conditions is satisfied$:$
\ben
\def\theenumi{(\Roman{enumi}$\mtho^{\rm B}$)}
\def\labelenumi{{\rm\theenumi}}
\itla{b1}\msur
$\mek$ is Borel linearizable -- moreover~\footnote
{\rm\ Harrington \ea\ \cite{hms} proved that any Borel \lo\ is 
Borel order isomorphic to a \lo\ 
$\stk{X}{\mel},$ where $X\sq 2^\al$ for some $\al<\omi$ and 
$\mel$ is the lexicographical order.}
in this case there are an ordinal\/ $\al<\omi$ and a Borel 
linearization\/ $h:\stk\cN\mek\lra\stk{2^\al}{\mel}\;;$

\itla{b2}
there exists a continuous \hop\/ $1-1$ map\/ 
$F:\stk{2^\om}\meo\lra\stk\cN\mek$ such that\/ 
$a\nEo b\imp F(a)\nmek F(b)$. 
\een
\ete

{\bfit Example 3.\/} Let 
$\wo=\ans{x\in\cN:x\,\hbox{ codes an ordinal}};$ for $x\in\wo$ 
let $|x|$ be the ordinal coded by $x.$ Define a $\is11$ \pqo\ 
$x\mek y$ iff either $y\not\in\wo$ or $x,\,y\in\wo$ and 
$|x|\<|y|.$ \ref{b1} is impossible for $\mek$ even via a non-Borel 
map $h$ since orders $\stk{2^\al}\mel,\msur$ $\al<\omi,$ do not 
admit strictly increasing \hbox{\dd\omi c}hains. \ref{b2} is also 
impossible via analytic maps $F$ by the restriction theorem. 
Thus Theorem~\ref{hb} fails for {\it analytic\/} 
relations.\vtm\hfill$\mtho\dashv$

Following ideas of Hjorth and Kechris \cite u, we involve longer 
orders, $2\lomi$ and $2^{\omi},$ to match the nature of analytic 
\pqo's. A set $A\sq 2\lomi$ will be 
called an {\it antichain\/} when it consists of pairwise 
\dd\sq incomparable elements. 

\bte
\label{ha}
Suppose that\/ ${\mek}$ is a\/ $\fs11$ \pqo\/ on\/ $\cN.$ Then at 
least one of the following two conditions is satisfied$:$
\ben
\def\theenumi{(\Roman{enumi}$\mtho^{\rm A}$)}
\def\labelenumi{{\rm\theenumi}}
\itla{1A}\msur
There is a linearization\/ 
$h:\stk\cN\mek\lra\stk{2^{\omi}}{\mel}$ such that for any\/ 
$\ga<\omi$ the map\/ $x\mapsto h(x)\res\ga$ is Borel, and has 
an\/ \dd\omi Borel~\footnote
{\rm\ Means: \dd\la Borel for an ordinal $\la<\omi$ which is not 
necessarily countable in $\kL[z]$.}
code in\/ $\kL[z]$ provided\/ $\mek$ is\/ 
$\is11[z].$ 
In addition in each of the two~\footnote
{\rm\ An obvious parallel with the ``Ulm classification'' theorem 
in Hjorth and Kechris~\cite{u} allows to conjecture that the 
additional assertion is also true in the assumption of the 
existence of ``sharps'', or an even weaker assumption in Friedman 
and Velickovic \cite{f}. However the most interesting problem is 
to prove the additional assertion in $\ZFC$.} 
following cases there is an antichain\/ $A\sq 2\lomi$ and a\/ 
$\fd12$ in the codes linearization\/ 
$h:\stk\cN\mek\lra\stk{A}{\mel}:$
\ben
\def\theenumii{(\alph{enumii})}
\def\labelenumii{{\rm\theenumii}}
\itla{1Abor} 
for any\/ $x$ the set\/ $[x]_\eee=\ans{y:y\eee x}$ is Borel$;$ 
\footnote
{\rm\ This applies \eg\ when $\mek$ is a \po. 
Recall that $x\eee y$ iff $x\mek y\cj y\mek x$.}

\itla{1Agen}
the universe is a set generic~\footnote
{\rm\ Via any kind of \underline{set} forcing. Compare with a 
theorem on thin $\fs11$ equivalence relations in 
Hjorth~\cite{h-thin}.
}
extension of a class\/ 
$\kL[z_0],\msur$ $z_0\in\cN$.
\een

\itla{2A}
As\/ {\rm\ref{b2}} of Theorem~\ref{hb}.
\een
\ete
Take notice that \ref{1A} and \ref{2A} here 
are compatible for instance in the assumption $\kV=\kL.$ There 
possibly exist reasonable sufficient conditions (like: all 
$\fd12$ sets are Lebesgue measurable) for \ref{1A} and \ref{2A} 
to be incompatible.

\punk*{\bfit Applications for analytic equivalence relarions}

Theorem \ref{ha} applies for analytic equivalence relations 
viewed as a particular case of \pqo's.

\bcor
\label e
Let\/ $\qE$ be a\/ $\fs11$ equivalence relation\/ on\/ $\cN.$ 
Then at least one of the following two conditions is satisfied$:$
\ben
\def\theenumi{(\Roman{enumi}$\mtho^{\rbox{E}}$)}
\def\labelenumi{{\rm\theenumi}}
\itla{1e}\msur
There is a map\/ $h:\cN\lra 2^{\omi}$ such that\/ 
${x\qE y}\eqv h(x)=h(y)$ and for any\/ 
$\ga<\omi$ the map\/ $x\mapsto h(x)\res\ga$ is Borel, and, 
provided\/ $\mek$ is\/ $\is11[z],$ has 
an\/ \hbox{\dd\omi B}orel code in\/ $\kL[z].$ 
In addition in each of the two following cases there is an 
antichain\/ $A\sq 2\lomi$ and a\/ $\fd12$ in the codes 
map\/ $h:\cN\lra A$ such that\/ ${x\qE y}\eqv h(x)=h(y):$
\ben
\def\theenumii{(\alph{enumii})}
\def\labelenumii{{\rm\theenumii}}
\itla{ebor} 
for any\/ $x$ the set\/ $[x]_{\qE}=\ans{y:y\qE x}$ is Borel$;$ 

\itla{egen}
the universe is a set generic extension of a class\/ 
$\kL[z_0],\msur$ $z_0\in\cN$.
\een

\itla{2e}
There exists a continuous\/ $1-1$ function\/ $F:2^\om\lra\cN$ such 
that\/ ${a\Eo b\eqv F(a)\qE F(b)}$.\qed
\een
\ecor
This result, with \ref{1e} in the additional form, has been 
obtained by Hjorth and Kechris~\cite{u} in the subcase \ref{ebor} 
(as well as in the assumption of existence of sharps), by 
Friedman and Velickovic~\cite{f} in a hypothesis connected with 
weakly compact cardinals, and by Kanovei~\cite{k-gen} in the 
subcase \ref{egen}. 

(Recall that a map $F$ as in \ref{2e} is called an 
{\it embedding\/} of $\Eo$ in $\qE$ -- a continuous embedding in 
this case. A map $h$ as in \ref{1e} is called a {\it reduction\/} 
of $\qE$ to the equality on $2^{\omi}$ or $2\lomi.$)

\bcor
\label{echar}
\footnote
{\rm\ Hjorth and Kechris told the author in April 1997 that they 
had known the result.}
Assume that\/ $\qE$ is a\/ $\is11[z]$ equivalence relation, 
$z\in\cN,$ and\/ {\rm\ref{2e}} of Corollary~\ref e fails. 
Then\/ $x\qE y$ iff we have\/ $x\in X\eqv y\in X$ for every\/ 
\dd\qE invariant Borel set\/ 
$X\sq \cN$ with an\/ \dd\omi Borel code in\/ $\kL[z]$.\qed
\ecor
Hjorth and Kechris~\cite u proved that any analytic $\qE$ which 
does not satisfy \ref{2e} of Corollary~\ref e admits an effective 
reduction $h:\cN\lra2^{\omi}$ (\ie\ we have $x\qE y\eqv h(x)=h(y)$), 
however it is not clear whether the property mentioned in 
\ref{1e} of Corollary~\ref e holds for the reduction given in 
\cite u and equally whether the reduction in \cite u directly 
leads to Corollary~\ref{echar}. 

\subsubsection*{\bfit Organization of the proofs}

The following theorem stands behind the results above. 

Recall that if $T$ is a tree on $\om\ti\om\ti\rdi$ then 
\dm
[T]=\ans{\ang{x,y,f}\in \cN^2\ti\la^\om:\kaz m\;
T(x\res m,y\res m,f\res m)}
\dm
and $\pr[T]=\ans{\ang{x,y}:\sus f\:[T](x,y,f)}$.

\bte
\label h
Let\/ $\om\<\rdi<\omi.$ 
Suppose that\/ $T$ and\/ $S$ are trees on\/ $\om\ti\om\ti\rdi$ 
such that the sets\/ ${\mekt}=\pr[T]\sq{\meks}=\dop\pr[S]$ are 
\pqo's on\/ $\cN.$ Then at least one of the following two 
conditions is satisfied$:$
\ben
\def\theenumi{(\Roman{enumi})}
\def\labelenumi{{\rm\theenumi}}
\itla1
There are\/ $\al<\omi$ and a\/ 
\dd\omi Borel coded in\/ $\kL[T,S]$ \hop\ 
map\/ $h:\stk\cN\mekt\lra\stk{2^\al}{\mel}$ 
such that\/ $h(x)=h(y)\imp x\ees y$.

\itla2
There is a continuous\/ $1-1$ \hop\ map\/ 
$F:\stk{2^\om}\meo\lra\stk\cN\mekt$ such that 
$a\nEo b\imp F(a)\nmeks F(b)$. 
\een
\ete

The principal technical scheme of the proofs goes back to the 
papers of Harrington and Shelah \cite{hs}, Shelah~\cite s, and 
Horth~\cite{h94} containing theorems on 
bi-$\mtho\kpa$-Souslin equivalence and order relations. 
However our version of the technique is free of any use of 
model theory including admissible sets. 

On the other hand we exploit several technical achievements made 
in the study of the Borel orders (Harrington~\ea~\cite{hms}, 
Louveau~\cite l) by means of the Gandi -- Harrington topology. 

The two technical schemes, the one we use and the one based on 
the Gandi -- Harrington topology, involve different kinds of 
``effective'' sets in the forcing, but 
have many common points in the construction of the proofs (like 
a similar definition of the ``regular'' and ``singular'' cases, a 
similar construction of splitting systems 
\etc), although differ in many details.

As a matter of fact the Gandi -- Harrington topology technique 
proves Theorem~\ref{hb} shorter than we do here (see 
Kanovei~\cite{k-b}), but it has problems with the analytic case as 
it does not capture the proper type of effectiveness. 

After some preliminaries in Section~\ref{prel} 
(including an effective version of the classical separation theorem) 
we introduce the dichotomy in 
Section~\ref{di}. Then the proof of Theorem~\ref h naturally 
develops itself in sections \ref{f}, \ref{pf}, \ref{E-b}, and 
\ref{abs} (where we show that \ref2 of Theorem~\ref h is 
Shoenfield--absolute).

Theorem~\ref{ha} (Section~\ref A) will require a reflection 
argument saying that an analytic \pqo\ has uncountably many indices 
for ``upper'' Borel approximations which are \pqo's, together with 
a delicate reasoning in the case of a generic universe, in 
Section~\ref S. 

\np

\parf{Preliminaries}
\label{prel}

The proof of Theorem~\ref h is the major part of this paper. 

We fix an ordinal $\rdi,\msur$ $\om\<\rdi<\omi,$ and   
trees $T,\,S\sq(\om\ti\om\ti\rdi)\lom.$ 

{\it Assume that both\/ $T$ and\/ $S$ are constructible\/}. 
\footnote
{\ Otherwise all entries of $\kL$ from now on have to be  
uniformly changed to $\kL[T,S]$.}

Suppose that ${\mekt}=\pr[T]\sq{\meks}=\dop\pr[S]$ are \pqo's on 
$\cN.$ Define $x\ees y$ iff $x\meks y\cj y\meks x$ and $x\eet y$ 
similarly.

\punk{Coding Borel sets}
\label C

We let $\lano$ be the infinitary language containing
\ben
\itemsep=1mm
\def\theenumi{(\roman{enumi})}
\def\labelenumi{\theenumi}
\itla{l3}
constant symbols $\dox,\,\doy,\,\doz,\,...$ for indefinite 
elements of $\cN=\om^\om$ and 
constant symbols $\dof,\,\dtg,\,...$ for indefinite 
elements of the set $\rdi^\om$;

\itla{l5}
elementary formulas of the form $\dox(k)=l$ and 
$\dof(k)=\al,$ where $k,\,l\in\om$ while $\al<\rdi$;

\itla{l6}
conjunctions and disjunctions of size $\<\rdi,$ 
together with the ordinary propositional connectives, but it is 
assumed that any formula contains only finitely many 
constant symbols mentioned in \ref{l3}.
\een
(Quantifiers are not allowed). 
Thus formulas in $\lano$ code \dd\rdii Borel subsets of 
spaces $\cN^m\ti(\rdi^\om)^n.$ 
For a formula, say, $\vpi(\dox,\dof)$ we put 
\dm
\mni\vpi=\ans{\ang{x,f}\in\cN\ti\rdi^\om:\vpi(x,f)}.
\dm 

\noi
For instance $[T](\dox,\doy,\dof)$ is a \dd\lano formula; we 
shall denote it by $\dox\mektf\doy.$ Similarly, the formula 
$\net [S](\dox,\doy,\dof)$ will be denoted by $\dox\meksf\doy.$ 
Formulas $\dox\neksf\doy,\msur$ $\dox\eesf\doy$ \etc\ are 
derivatives. Then 
\dm
x\mekt y\eqv\sus f\in\la^\om\:x\mektF y
\hspace{4mm}\hbox{and}\hspace{4mm}
x\meks y\eqv\kaz f\in\la^\om\:x\meksF y\,.
\dm


\punk{Consistency and separation}
\label{co}
\label{sep}

A formula $\vpi$ is {\it consistent\/} if it has a model, \ie\ 
becomes true after one suitably substitutes its constants by 
elements of $\cN$ and $\rdi^\om.$ 
 
A {\it theory\/} in $\lano$ will be any set of formulas of 
$\lano$ containing a common (finite) list\/ $S$ of constants of 
type \ref{l3}. (We shall usually 
consider {\it constructible\/} theories $\Phi\sq\lano$.) 
A theory is {\it consistent\/} if it has a model. A theory $\Phi$ 
is \dl{\it consistent\/} if every constructible subtheory 
$\Phi'\sq \Phi$ of cardinality $\<\rdi$ in $\kL$ is consistent. 

A theory $\Phi$ \dl{\it implies\/} a formula $\psi$ 
if $\Phi\cup\ans{\psi}$ is \dl inconsistent. Other statements like 
this are to be understood accordingly. 

The following theorem has a semblanse of the Craig interpolation 
theorem, but essentially it belongs to the type of 
{\it separation\/} theorems.

\bte
\label{tsep}
Suppose that\/ $\Phi(\dox,\doy,\dof,...)$ and 
$\Psi(\dox,\doy',\dof',...)$ are construct\-ible theories in\/ 
$\lano$ having\/ $\dox$ as the only common constant in the 
(finite) lists of constants. Assume that\/ 
$\Phi(\dox,\doy,\dof,...)\cup\Psi(\dox,\doy',\dof',...)$ 
is\/ \dl in\-consistent. Then there is a\/ \dd\lano formula\/ 
$\pi(\dox)$ \dl separating\/ $\Phi$ from\/ $\Psi$ in the sense 
that\/ $\Phi(\dox,...)$ \dl implies\/ $\pi(\dox)$ while\/ 
$\Psi(\dox,...)$ \dl implies\/ $\net\pi(\dox)$.
\ete
\bpf
First of all we can assume that $\Phi$ and 
$\Psi$ consist of single formulas, resp.\ 
$\vpi(\dox,\doy,\dof,...)$ and $\psi(\dox,\doy',\dof',...).$ Let, 
for the sake of simplicity, $\vpi$ be $\vpi(\dox,\doy)$ and 
$\psi$ be $\psi(\dox,\dof).$ Consider the sets
\dm
P=\mni{\vpi}=\ans{\ang{x,y}:\vpi(x,y)}
\hspace{1mm},\hspace{4mm}
Q=\mni{\psi}=\ans{\ang{x,f}\in\cN\ti\rdi^\om:\psi(x,f)}\,.
\dm
The projections $X=\ans{x:\sus y\:P(x,y)},\msur$ 
$Y=\ans{x:\sus f\:Q(x,f)}$ are disjoint $\fs11$ sets by the 
inconsistency assumption,  
hence by the classical separation theorem they can be separated 
by a Borel set. Moreover as we demonstrated in 
Kanovei~\cite{k-g} (Theorem 7) in this case the separating set 
can be defined in the form $B=\mni\pi$ for an 
appropriate \dd\lano formula $\pi(\dox)$.
\epf

\punk{Hulls}
\label{flas}

By $\kaf(\dox)$ we shall denote the (constructible) collection of 
all \dd\lano formulas $\vpi(\dox)\in\kL.$ 
For a theory $\Phi(\dox,\doy,...),$ $\hax{\Phi(\dox,\doy,...)}$ 
will be the set of all formulas $\vpi(\dox)\in\kaf(\dox)$ which 
are \dl implied by $\Phi(\dox,\doy,...).$ 

\ble
\label{str-o}
Suppose that\/ $\Pi(\dox,...)$ is a constructible theory in\/ 
$\lano$ while $\Rho(\dox)\sq\kaf(\dox)$ is also constructible 
and satisfies\/ $\hax{\Pi(\dox,...)}\sq\Rho(\dox).$ Then the 
theory\/ $\Pi'(\dox,...)=\Pi(\dox,...)\cup \Rho(\dox)$ 
satisfies\/ $\hax{\Pi'(\dox,...)}=\hax{\Rho(\dox)}$. 
\ele
\bpf
Prove that $\hax{\Pi'(\dox,...)}\sq\hax{\Rho(\dox)}$ (the 
nontrivial direction). Let 
$\psi(\dox)\in\hax{\Pi'(\dox,...)}.$ By definition there is 
a set $\Psi(\dox)\in\kL,\msur$ $\Psi(\dox)\sq\Rho(\dox),$ 
of cardinality $\<\rdi$ in $\kL,$ such that 
$\Pi(\dox,...)\cup\Psi(\dox)$ \dl implies $\psi(\dox).$ We 
conclude that the formula $(\bwe\Psi(\dox))\imp\psi(\dox)$ 
belongs to $\hax{\Pi(\dox,...)},$ hence to $\Rho(\dox),$ 
which guarantees $\psi(\dox)\in\hax{\Rho(\dox)}$. 
\epf

\np

\parf{Monotone Borel functions and the dichotomy}
\label{di}

To introduce the dichotomy we have to extend the language 
$\lano$ by Borel functions mapping $\cN$ in a set of the form 
$2^\al,$ where $\al<(\rdi^+)^\kL.$ If such an $\al$ is fixed, 
let a {\it function code\/} be a sequence of the form 
$\vvi=\ang{\vpi_\ga(\dox):\ga<\al}$ where each $\vpi_\ga$ is a 
\dd\lano formula. Such a sequence defines a function 
$h_\vvi:\cN\lra 2^\al$ so that  
$h_\vvi(x)(\ga)=1$ iff $\vpi_\ga(x).$ 

Define $\rH_\al$ to be the set of all \hop\ maps  
$h_\vvi:\stk\cN\mekt\lra\stk{2^\al}{\mel},$ where $\vvi$ is a 
{\it constructible\/} function code. 
%
Define $\rH=\bigcup_{\al<(\rdi^+)^\kL}\rH_\al.$ Thus every 
function in $\rH$ is an \dd\omi Borel (even \dd{\la\pone}Borel), 
coded in $\kL,$ map from $\cN$ to some $2^\al,\msur$ 
$\al<(\rdi^+)^\kL,$ satisfying ${x\mekt y}\imp {h(x)\mel h(y)}$.

But as a matter of fact functions in $\rH$ will be used only 
via equalities of the form $h(x)=h(y)$ where $h=h_\vvi\in\rH_\al$ 
for some $\al<(\rdi^+)^\kL,$ viewed as shorthand for 
$\bwe_{\ga<\al}(\vpi_\ga(x)\eqv\vpi_\ga(y))$.

Let $\dox\eqh\doy$ be the theory $\ans{h(\dox)=h(\doy):h\in\rH}.$ 
(Thus $\eqh$ defines an equivalence relation extending $\eet$.) 

We have two cases.~\footnote
{\label{diss}\ There is a point of dissatisfaction in the 
distribution on the two cases we use. It would be more natural 
to define Case~1 as that $\dox\meh\doy$ \dl implies 
$\dox\meks\doy,$ where $\dox\meh\doy$ is the theory 
$\ans{h(\dox)\mel h(\doy):h\in\rH},$ which would 
improve \ref1 of Theorem~\ref{h} to the existence of a \hop\ 
map satisfying $h(x)\leq h(y)\imp x\meks y$. However then 
the arguments for Case~2, especially the key lemmas in the next 
section, do not go through.}
\vtm

{\bfit Case 1\/}{\bf:} the theory $\dox\eqh \doy$ \dl implies 
$\dox\eesf\doy$.\vtm 

Then clearly there is a single function $h\in \rH$ such that 
$h(\dox)=h(\doy)$ already implies $\dox\eesf\doy.$
Then $h$ satisfies \ref1 of Theorem~\ref h. 
\vtm

{\bfit Case 2\/}{\bf:} the theory 
$(\dox\eqh\doy)\cup\ans{\dox\nesf\doy}$ is \dl consistent.\vtm  

Assuming this we shall work towards \ref2 of Theorem~\ref h. We 
begin with a study of an important class of ``conditionally 
downward closed'' formulas.

Let $\kah(\dox)$ be the (constructible as above) set of all 
formulas $\eta(\dox)\in\kaf(\dox)$ satisfying the following: 
$\dox\eqh\doxp$ \dl implies 
$\eta(\dox)\cj\doxp\mektf\dox\imp\eta(\doxp)$. 

\ble
\label{dkle}
Suppose that\/ $\eta(\dox)\in\kah(\dox).$ Then there is a 
function\/ $h\in \rH_{\al+1}$ for some\/ $\al<(\rdi^+)^\kL$ such 
that\/ $\eta(x)\eqv h(x)=0.$ In particular the 
theory\/ $\eta(\dox)\cj \net\eta(\doy)\cj \dox\eqh\doy$ is then 
\dl inconsistent.
\ele
\bpf
By definition there exists a function 
$h_0\in H_\al$ for some $\al<(\rdi^+)^\kL$ satisfying 
$h_0(x)=h_0(y)\cj\eta(x)\cj x'\mekt x\imp \eta(x').$ 
Define $h(x)=h_0(x)\we 0$ whenever $\eta(x)$ and 
$h(x)=h_0(x)\we 1$ otherwise. 
\epf

For a theory $\Phi(\dox,\doy,...),$ let 
$\hhx{\Phi(\dox,\doy,...)}=\hax{\Phi(\dox,\doy,...)}\cap
\kah(\dox)$.

\np

\parf{The basic forcing}
\label{f}

Let $\Xi(\dox)$ denote the (constructible) set of all formulas 
$\xi(\dox)\in\kaf(\dox)$ which are \dl implied by the theory 
$\dox\eqh\doy\cup\ans{\dox\nesf\doy}$. 

Let $\dP$ be the set of all \dl consistent theories 
$\Pi(\dox)\in\kL,\msur$ $\Pi\sq\kaf(\dox),$ including 
$\Xi(\dox).$ Then $\dP\in\kL,$ so we can view it as a forcing 
notion over $\kL$.

\ble
\label{dpa-b}
Let\/ $G\sq\dP$ be\/ \dd\dP generic over\/ $\kL.$ Then there is a 
unique real\/ $x=x_G\in\kL[G]$ such that\/ $\pi(x)$ holds in\/ 
$\kL[G]$ for any formula\/ $\pi(\dox)\in \bigcup G$. 
\ele
\bpf
Note that, for any $n,$ the set $D_n$ of all theories 
$\Pi(\dox)\in \dP$ which include $\dox(n)=m$ for some $m$ is 
dense in $\dP$ and belongs to $\kL,$ hence 
$D_n\cap G\not=\emps.$ The rest of the proof is standard.
\epf

\punk{Key lemmas}
\label{key}

\ble
\label{k1}
Let\/ $\Pi(\dox)$ be a theory in\/ $\dP.$ Then the theory\/
$\Phi_\Pi(\dox,\doy,\dof)\eqdf\Pi(\dox)\cup\Pi(\doy)\cup
\dox\eqh\doy\cup\ans{\dox\neksf\doy}$ 
is \dl consistent and satisfies the equalities\/ 
$\hax{\Phi_\Pi(\dox,\doy,\dof)}=\hax{\Pi(\dox)}$ and\/ 
$\hay{\Phi_\Pi(\dox,\doy,\dof)}=\hay{\Pi(\doy)}$. 
\ele
\bpf
Let us first prove the consistency. Otherwise there is  
a formula $\pi(\dox)\in\hax{\Pi(\dox)}$ and a function $h\in\rH$ 
such that the formula 
$\pi(\dox)\cj\pi(\doy)\cj h(\dox)=h(\doy)\cj \dox\neksf\doy$ 
is \dl inconsistent. 

The plan is to find functions $h',\,h''\in\rH$ such that the 
theories $\pi(\dox)\cj h'(\dox)=h'(\doy)\cj \doy\neksf\dox$ and 
$\pi(\dox)\cj h''(\dox)=h''(\doy)\cj \dox\neksf\doy$
are \dl in\-consistent: then the formula $\net\pi(\dox)$ 
belongs to $\Xi,$ which is a contradiction because $\Pi$ 
includes $\Xi$. 

Consider the first theory (the other one is similar). 
By Lemma~\ref{dkle} it 
suffices to get a formula $\psi(\dox)\in\kah(\dox)$ such that 
$X=\mni{\pi}\sq U=\mni\psi$ and, for all 
$x\in X$ and $u\in U,$ $h(x)=h(u)$ implies $u\meks x$.  

Let $Z=\ans{z:\kaz x\in X\:(h(z)=h(x)\imp z\meks x)}.$ Then  
$X\sq Z$ by the inconsistency assumption above. 

Define a sequence of sets 
$X=X_0\sq U_0\sq X_1\sq U_1\sq \dots\sq Z$ and formulas 
$\vpi_n(\dox)\in\kaf$ so that 
$U_n=\ans{u:\sus x\in X_n\:(h(x)=h(u)\cj u\mekt x)},$ 
$X_n=\mni{\vpi_n(\dox)},$ and the sequence of formulas $\vpi_n$ 
is constructible.

Now, $\psi(\dox)=\bve_n\vpi_n(\dox)$ is the formula required.  
(Note that $\mni\psi=\bigcup_n X_n=\bigcup_n U_n$.)  
It remains to carry out the construction of $X_n,\,U_n,\,\vpi_n$. 

Suppose that $X_n=\mni{\vpi_n(\dox)}\sq Z$ has been defined. 
Define $U_n$ by the equality above. Then $X_n\sq U_n,$ and 
$U_n\sq Z.$ 
(Assume that $u\in U_n,$ so $u\mekt x$ for some $x\in X_n$ 
satisfying $h(x)=h(u).$ Take any $x'\in X$ such that $h(x')=h(u)$ 
and prove $u\meks x'.$ First of all $h(x)=h(x')$ hence $x\meks x'$ 
because $x\in X_n\sq Z.$ Now $u\meks x'$ as $u\mekt x$.) 

Theorem~\ref{tsep} yields a formula 
$\ba(\dox)\in\kaf$ such that the set $B=\mni\ba$ 
satisfies $U_n\sq B\sq Z.$ Take $X_{n+1}=B$ and $\vpi_{n+1}=\ba$. 

Finally, as the choice of the formulas $\vpi_n$ can be forced 
in $\kL,$ the sequence of formulas can be chosen in $\kL.$ 
This ends the proof of the consistency. 

The equality 
$\hax{\Phi_\Pi(\dox,\doy,\dof)}=\hax{\Pi(\dox)}$ does not cause 
much trouble. Indeed suppose that 
$\pi(\dox)\in\kaf(\dox)\setminus\hax{\Pi(\dox)}$ (the nontrivial 
direction). Then the theory 
$\Pi'(\dox)=\Pi(\dox)\cup\ans{\net\pi(\dox)}$ 
is \dl consistent, hence belongs to $\dP.$ It follows from the 
above that 
$\Phi_\Pi(\dox,\doy,\dof)\cup\ans{\net\pi(\dox)}$ 
is \dl consistent as well, hence 
$\pi(\dox)\not\in\hax{\Phi_\Pi(\dox,\doy,\dof)}$. 
\epf

\ble
\label{k2}
Assume that\/ $\Pi(\dox)$ and\/ $\Rho(\dox)$ belong to\/ 
$\dP,$ and\/ $\hhx{\Pi(\dox)}=\hhx{\Rho(\dox)}.$ Then\/ 
$\Psi_{\Pi\Rho}(\dox,\doy,\dof)\eqdf\Pi(\dox)\cup \Rho(\doy)\cup 
\dox\eqh\doy\cup \ans{\dox\mektf\doy}$ 
is a \dl consistent theory satisfying the equalities\/ 
$\hax{\Psi_{\Pi\Rho}(\dox,\doy,\dof)}=\hax{\Pi(\dox)}$ and\/ 
$\hay{\Psi_{\Pi\Rho}(\dox,\doy,\dof)}=\hay{\Rho(\doy)}$. 
\ele
\bpf
As in the previous lemma, it suffices to prove the consistency. 
Suppose otherwise. Then there exist formulas 
$\pi(\dox)\in\hax{\Pi(\dox)}$ and $\rho(\doy)\in\hay{\Rho(\doy)},$ 
and a function $h\in\rH$ such that the formula 
$\pi(\dox)\cj\rho(\doy)\cj h(\dox)=h(\doy)\cj \dox\mektf\doy$ 
is inconsistent. In other words we have $x\nmekt y$ whenever 
$x\in X=\mni\pi$ and $y\in Y=\mni\rho$ satisfy $h(x)=h(y)$.

Define 
$Z=\ans{z:\kaz y\in Y\:(h(y)=h(z)\imp z\nmekt y)},$ so that 
$X\sq Z$ but $Y\cap Z=\emps.$ The same iterated procedure as in 
the proof of Lemma~\ref{k1} (with 
$U_n=\ans{u:\sus x\in X_n\:(h(x)=h(u)\cj x\mekt u)}$ 
yields a formula $\psi(\dox)\in\kah(\dox)$ 
such that the set $U=\dop{\mni\psi}$ satisfies $X\sq U\sq Z.$ 
But this contradicts the assumption 
$\hhx{\Pi(\dox)}=\hhx{\Rho(\dox)}$.
\epf

\bcor
\label{k3}
Suppose that\/ $\Pi(\dox),\msur$ $\Rho(\dox)$ belong to\/ 
$\dP$ and\/ $\hhx{\Pi(\dox)}=\hhx{\Rho(\dox)}.$ Then\/
$\rU(\dox,\doy,\dof)\eqdf\Pi(\dox)\cup\Rho(\doy)\cup
\dox\eqh\doy\cup\ans{\dox\neksf\doy}$ 
is a \dl consistent theory satisfying the equalities\/ 
$\hax{\rU(\dox,\doy,\dof)}=\hax{\Pi(\dox)}$ and\/ 
$\hay{\rU(\dox,\doy,\dof)}=\hay{\Rho(\doy)}$. 
\ecor
\bpf
It suffices, as above, to prove the consistency. 
Suppose otherwise. Then the theory 
$\Phi_\Pi(\dox,\doz,\dof)\cup\Psi_{\Rho\Pi}(\doy,\doz,\dtg)$ 
is \dl inconsistent as well. (Otherwise we have 
reals $x,\,y,\,z$ satisfying $\Pi(x),\msur$ $\Pi(z),\msur$ 
$x\eqh z,\msur$ $\Rho(y),$ and $y\eqh z,$ hence $x\eqh z,$ and, 
in addition, $x\nmeks z$ and $y\mekt z,$ hence $x\nmeks y$.)
Theorem~\ref{tsep} yields a formula $\pi(\doz)\in\kaf(\doz)$ 
\dl implied by $\Phi_\Pi(\dox,\doz,\dof)$ but \dl inconsistent 
with $\Psi_{\Rho\Pi}(\doy,\doz,\dtg),$ which is a contradiction 
as 
$\haz{\Phi_\Pi(\dox,\doz,\dof)}=\haz{\Pi(\doz)}
=\haz{\Psi_{\Rho\Pi}(\doy,\doz,\dtg)}$ 
by lemmas ~\ref{k1}, \ref{k2}. 
\epf

\bcor
\label{kd}
Suppose that\/ $\Pi(\dox),\msur$ $\Rho(\dox)$ belong to\/ 
$\dP$ and\/ $\hhx{\Pi(\dox)}=\hhx{\Rho(\dox)}.$ Then there are 
theories\/ $\Pi'(\dox),\,\Rho'(\dox)\in \dP$ such that\/ 
$\Pi\sq\Pi',\msur$ $\Rho\sq\Rho',$ still\/ 
$\hhx{\Pi'(\dox)}=\hhx{\Rho'(\dox)},$ and\/ 
$\Pi'(\dox)\cup\Rho'(\dox)$ is\/ \dl inconsistent. 
\ecor
\bpf
The theory 
$\Phi(\dox,\doy)\eqdf\Pi(\dox)\cup\Rho(\dox)\cup(\dox\eqh\doy)
\cup\ans{\dox\not=\doy}$ 
is \dl consistent by the previous corollary. It easily follows 
that there exist numbers $m$ and $k_x\not=k_y$ such that 
$\Phi'(\dox,\doy)\eqdf
\Phi(\dox,\doy)\cup\ans{\dox(m)=k_x}\cup\ans{\doy(m)=k_y}$ 
is still \dl consistent. Now set 
$\Pi'(\dox)=\hhx{\Phi'(\dox,\doy)}$ and 
$\Rho'(\doy)=\hhy{\Phi'(\dox,\doy)}$.
\epf

\punk{Two-dimentional modifications}
\label{2dim}

There are two related forcing notions which produce pairs 
of reals. 

Let $\dpw$ be the (constructible) set of all 
\dl consistent theories $\Da(\dox,\doy)\in\kL$ such that 
$\Xi(\dox)\cup\Xi(\doy)\sq \Da(\dox,\doy)$. \vtm

\noi
{\bfit First modification.\/} 
Recall that ${\mekt}=\pr[T].$ The idea is to define 
a forcing which leads to pairs of reals satisfying $x\mekt y.$ 

We let $\dpwt(\dox,\doy,\dof)$ be the set of all \dl consistent 
theories $\Tau(\dox,\doy,\dof)$ of the form 
${\Da(\dox,\doy) \cup \rF \cup \dox\eqh\doy\cup
\ans{\dox\mektf\doy}},$ 
where $\Da\in\dpw$ and $\rF$ is a finite collection of formulas 
$\dof(k)=\al$ (where $k\in\om$ and $\al<\rdi$).

For instance the theory 
$\Xi(\dox)\cup\Xi(\doy)\cup\dox\eqh\doy\cup\ans{\dox\mektf\doy}$ 
(which is \dl consistent by Lem\-ma~\ref{k2}) belongs to $\dpwt$. 

\ble
\label{dpat-b}
Let\/ $G\sq\dpwt$ be\/ \dd\dpwt generic over\/ $\kL.$ There is 
a unique triple\/ 
$\ang{x,y,f}\in\kL[G]\cap(\cN\ti\cN\ti\rdi^\om)$ such that\/ 
$\tau(x,y,f)$ holds for any formula\/ 
$\tau(\dox,\doy,\dof)\in \bigcup G.$ 
In particular we have\/ $x\mekt y$. 
\ele
\bpf
Analogous to Lemma~\ref{dpa-b}.
\epf

\noi
{\bfit Second modification.\/} 
Recall that ${\meks}=\dop\pr[S].$ Now the intension is to force 
pairs of reals $\ang{x,y}$ satisfying 
$x\nmeks y.$ We let $\dpws(\dox,\doy,\dof)$ be the set of all 
\dl consistent theories $\Sg(\dox,\doy,\dof)$ of the form 
$
\Da(\dox,\doy) \cup \rF \cup \dox\eqh\doy \cup 
\ans{\dox\neksf\doy},
$ 
where $\rF$ and $\Da$ are as in the definition of $\dpwt$.

To see that $\dpws\not=\emps$ note that   
$\Xi(\dox)\cup\Xi(\doy)\cup\dox\eqh\doy\cup\ans{\dox\neksf\doy}$ 
is a \dl consistent theory by Lemma~\ref{k1}.

\ble
\label{dpas-b}
Let\/ $G\sq\dpws$ be\/ \dd\dpws generic over\/ $\kL.$ There is 
a unique triple\/ $\ang{x,y,f}\in\kL[G]\cap(\cN\ti\cN\ti\rdi^\om)$ 
such that\/ $\sg(x,y,f)$ holds for any elementary formula\/  
$\sg(\dox,\doy,\dof)\in \bigcup G.$ 
In particular we have\/ $x\nmeks y$.
\ele
\bpf
Analogous to Lemma~\ref{dpa-b}. 
\epf

\np

\parf{The product forcing}
\label{pf}

The forcing notion $\dpws$ executes too a tight control over 
generic reals. Fortunately generic \dd\meks incomparable 
pairs can be obtained by another forcing, which connects the 
components in a much looser way, so that it is rather a kind 
of product forcing, with the factors equal to $\dP$.

We let $\dpwo$ be the set of all theories $\Ups(\dox,\doz)$ 
of the form $\Pi(\dox)\cup\Rho(\doz),$ 
where $\Pi$ and $\Rho$ belong to $\dP$ and satisfy 
$\hhx{\Pi(\dox)}=\hhx{\Rho(\dox)}$. The set $\dpwo$ 
is constructible and non-empty.

\bte
\label{dpao-b}
Let\/ $G\sq\dpwo$ be\/ \dd\dpwo generic over\/ $\kV.$ There is 
a unique pair\/ ${\ang{x,z}\in\kV[G]\cap\cN^2}$ such that\/ 
$\ups(x,z)$ holds for any formula\/ $\ups(\dox,\doz)\in\bigcup G.$ 
Moreover we have\/ $x\nmeks z$.
\ete

Pairs $\ang{x,z}$ as in the theorem will be denoted by 
$\ang{x_G,z_G}$ and called \dd\dpwo{\it generic over\/ $\kV$}. 
($\kV$ is the universe of all sets as usual.)

\bpf
Let us concentrate on the proof that $x_G$ and $z_G$ are 
\dd\meks incompar\-able; the rest is analogous to the above. 

Suppose on the contrary that $x_G\meks z_G$ is \dd{\dpwo}forced 
over $\kV$ by a ``condition '' 
$\Ups_0(\dox,\doz)=\Pi_0(\dox)\cup\Rho_0(\doz)\in\dpwo,$ where 
$\Pi_0$ and $\Rho_0$ belong to $\dpw$ and 
satisfy $\hhx{\Pi_0(\dox)}=\hhx{\Rho_0(\dox)}$. 

We shall define a generic ``rectangle'' of reals 
$x,\,z,\,x',\,z',$ such that the following will be forced: 
$x\meks z$ and $x'\meks z'$ -- by the forcing assumption, 
$z\mekt x'$ -- by Lemma~\ref{dpat-b}, and $x\nmeks z'$  -- 
by Lemma~\ref{dpas-b}, which is a contradiction. 
The forcing $\gP$ used to get a required ``rectangle'' consists 
of forcing conditions of the following general form:
\dm
\gp=
\ang{\Ups(\dox,\doz),\Tau(\doz,\doxp,\dof),\Ups'(\doxp,\dozp),
\Sg(\dox,\dozp,\dof)}
\dm 
such that the theories $\Ups(\dox,\doz)=\Pi(\dox)\cup\Rho(\doz)$ 
and $\Ups'(\doxp,\dozp)=\Pi'(\doxp)\cup\Rho(\dozp)$ belong to 
$\dpwo,$ $\Tau$ belongs to $\dpwt,$ $\Sg$ belongs 
to $\dpws,$ and 
\dm
\bay{rclcrcl}
\Pi(\dox)&=&\hax{\Sg(\dox,\dozp,\dof}
& \hbox{, }\, &
\Rho(\doz)&=&\haz{\Tau(\doz,\doxp,\dof)} \\[2mm]

\Pi'(\doxp)&=&\haxp{\Tau(\doz,\doxp,\dof)} 
& \hbox{, }\, &
\Rho'(\dozp)&=&\hazp{\Sg(\dox,\dozp,\dof)} 
\eay
\dm
Order $\gP$ componentwise: $\gp_1$ is stronger than $\gp_2$ 
iff $\Ups_2(\dox,\doz)\sq \Ups_1(\dox,\doz),\msur$ 
$\Tau_2(\doz,\doxp,\dof)\sq \Tau_1(\doz,\doxp,\dof),\msur$ 
$\Ups'_2\sq \Ups'_1,$ and $\Sg_2\sq\Sg_1$.

To get a condition in $\gP$ we start with the given theory 
$\Ups_0(\dox,\doz)=\Pi_0(\dox)\cup\Rho_0(\doz)\in\dpwo.$ 
By definition $\hhx{\Pi_0(\dox)}=\hhx{\Rho_0(\dox)}$. 

It can be supposed that 
$\hax{\Pi_0(\dox)}=\Pi_0(\dox)$ (otherwise replace $\Pi_0(\dox)$ 
by $\hax{\Pi_0(\dox)}$) and $\haz{\Rho_0(\doz)}=\Rho_0(\doz).$ 

The theory 
$\Sg_0(\dox,\dozp,\dof)\eqdf{\Pi_0(\dox)\cup\Rho_0(\dozp)\cup
{\dox\eqh\dozp}\cup {\dox\neksf\dozp}}$ 
belongs to $\dpws(\dox,\dozp)$ and satisfies the equalities 
$\hax{\Sg_0(\dox,\dozp,\dof)}=\Pi_0(\dox)$ and 
$\hazp{\Sg_0(\dox,\dozp,\dof)}=\Rho_0(\dozp)$ by 
Corollary~\ref{k3}. 
Similarly, by Lemma~\ref{k2}, 
the theory 
$\Tau_0(\doz,\doxp,\dof)\eqdf{\Rho_0(\doz)\cup\Pi_0(\doxp)\cup
{\doz\eqh\doxp}\cup{\doz\mektf\doxp}}$ 
belongs to $\dpwt$ and satisfies 
$\haz{\Tau_0(\doz,\doxp,\dof)}=\Rho_0(\doz)$ and 
$\haxp{\Tau_0(\doz,\doxp,\dof)}=\Pi_0(\doxp)$. 

Now $\gp_0=\ang{\Ups_0(\dox,\doz),\Tau_0(\doz,\doxp,\dof),
\Ups_0(\doxp,\dozp), \Sg_0(\dox,\dozp,\dof)}$
belongs to $\gP$.

\bass 
Suppose that\/ $\gp=\ang{\Ups,\Tau,\Ups',\Sg}\in\gP,\msur$  
$\Ups_1(\dox,\doz)\in\dpwo,$ and\/ 
$\Ups(\dox,\doz)\sq\Ups_1(\dox,\doz).$ Then there exists\/ 
$\gp_1=\ang{\Ups_1,\Tau_1,\Ups'_1,\Sg_1}\in\gP$ stronger 
than\/ $\gp$ {\rm(}\ie\ $\Tau\sq\Tau_1,\msur$ 
$\Ups'\sq\Ups'_1,$ and\/ $\Sg\sq\Sg_1${\rm).} 

The same is true when we strengthen any of the other three
components.
\eass
\bpfi
By definition $\Ups_1(\dox,\doz)=\Pi_1(\dox)\cup\Rho_1(\doz)$ 
where $\Pi_1$ and $\Rho_1$ belong to $\dpw$ and 
$\hhx{\Pi_1(\dox)}=\hhx{\Rho_1(\dox)}$. 

Let  
$\Tau_1(\doz,\doxp,\dof)=\Tau(\doz,\doxp,\dof)\cup\Rho_1(\doz).$ 
Lemma~\ref{str-o} yields 
$\hhz{\Tau_1(\doz,\doxp,\dof)}=\hhz{\Rho_1(\doz)}.$ 
Let $\Pi'_1(\doxp)=\haxp{\Tau_1(\doz,\doxp,\dof)}.$ Now 
$\hhz{\Pi'_1(\doz)}=\hhz{\Rho_1(\doz)}$ by 
Lemma~\ref{dkle} because even $\Tau(\doz,\doxp,\dof)$ includes 
$\doz\eqh\doxp$ by the 
definition of $\dpwt.$ Hence 
$\hhx{\Pi'_1(\dox)}=\hhx{\Pi_1(\dox)}$ by the above.  

We define $\Rho'_1(\dozp)$ using the other ark of the rectangle. 
Let 
$\Sg_1(\dox,\dozp,\dof)=\Sg(\dox,\dozp,\dof)\cup\Pi_1(\dox)$ 
and $\Rho'_1(\dozp)=\hazp{\Sg_1(\dox,\dozp,\dof)}.$ 
Then $\hhx{\Rho'_1(\dox)}=\hhx{\Pi_1(\dox)}$ by lemmas \ref{str-o} 
and \ref{dkle} as above. Thus in particular we have 
$\hhxp{\Rho'_1(\doxp)}=\hhxp{\Pi'_1(\doxp)},$ so that 
$\Ups'_1(\doxp,\dozp)=\Pi'_1(\doxp)\cup \Rho'_1(\dozp)$ 
belongs to $\dpwo,$ closing the diagram. Clearly 
$\Ups\sq\Ups_1.$ It easily follows from the construction that 
$\gp_1=\ang{\Ups_1,\Tau_1,\Ups'_1,\Sg_1}\in\gP$ is as 
required. 
\epfi

We continue the proof of Theorem~\ref{dpao-b}. Consider a 
\dd\gP generic extension $\kV[\gG]$ such that the generic set 
$\gG\sq\gP$ contains the condition $\gp_0$ defined above. It 
easily follows from the assertion just proved that $\gG$ results 
in a ``rectangle'' of reals $x,\,z,\,x',\,z'\in \kV[\gG]$ such 
that 

\ben
\def\theenumi{\arabic{enumi})}
\def\labelenumi{\theenumi}
\itla{zx}
The pair $\ang{z,x'}$ is \dd\dpwt generic over 
$\kV,$ therefore we have $z\mekt x'$ in $\kV[\gG]$ 
by Lemma~\ref{dpat-b}.

\itla{xy}
The pair $\ang{x,z'}$ is \dd\dpws generic over $\kV,$ therefore we 
have $x\nmeks y$ in $\kV[\gG]$ by Lemma~\ref{dpas-b}.

\itla{xz} 
The pairs $\ang{x,z}$ and $\ang{x',z'}$ are \dd\dpwo generic over 
$\kV,$ and moreover, the corresponding generic subsets of $\dpwo$ 
contain the ``condition'' $\Ups_0(\dox,\doz)$ fixed above; hence 
we have $x\meks z$ and $x'\meks z'$ in $\kV[\gG]$ by the choice of 
$\Ups_0$.
\een
This is a contradiction because ${\mekt}\sq{\meks}$ in $\kV[\gG]$ 
by absoluteness.
\epf

\np

\parf{The construction of an embedding}
\label{E-b}

We are going to define a continuous $1-1$ map $F:2^\om\lra\cN$ 
satisfying \ref2 of Theorem~\ref h. Our strategy will be to prove 
the existence of such a map in a \dd\kpa collapse generic extension 
$\kvp$ of $\kV,$ the universe of all sets, where $\kpa$ is the 
cardinal $2^{2^{\aleph_0}}$ in $\kV.$ This suffices to conclude 
that \ref2 of Theorem~\ref h holds in $\kV$ by Lemma~\ref{l-abs} 
of the next section.

\punk{Generic splitting family of theories}
\label{sf}

Let a {\it crucial pair\/} be any (ordered) pair $\ang{u,v}$ 
such that $u,\,v\in 2^m$ for some $m$ and 
$u=1^k\we 0\we w,\msur$ $v=0^k\we1\we w$ for some $k<m$ 
and $w\in 2^{m-k-1}.$ 

By the choice of $\kvp$ the sets $\dP,\msur$ 
$\dpwt,\msur$ $\dpwo$ have only \dd\kvp countably many subsets 
in $\kV.$ 
Let $\ans{D(n):n\in\om},\msur$ $\ans{D^2(n):n\in\om},\msur$ 
$\ans{D_2(n):n\in\om}$ be enumerations, in $\kvp,$ of the 
collections of all dense (by {\it dense\/} we mean 
{\it open dense\/}) subsets of resp.\ $\dP,\msur$ $\dpwo,\msur$ 
$\dpwt.$ It will be assumed that each dence set has infinitely 
many indices in the relevant enumeration. 


We shall define, in $\kvp,$ a family of theories 
$\Pi_u(\dox)\in\dP$ (where $u\in 2\lom$) and 
$\Tau_{uv}(\dox,\doy,\dof)\in\dpwt$ 
(where $\ang{u,v}$ is a crucial pair in some $2^n$) satisfying 
the following conditions, for all $u\in 2\lom$ and $i=0,1$:

\ben
\def\theenumi{(\roman{enumi})}
\def\labelenumi{\theenumi}
\itla 0 
$\Pi_u\in D(n)$ whenever $u\in 2^n$; \ 
$\msur\Pi_u(\dox)\sq\Pi_{u\we i}(\dox)$;

\itla o
if $\ang{u,v}$ is a crucial pair in $2^n$ then 
$\Tau_{uv}\in D_2(n)$; \ $\msur\Tau_{uv}(\dox,\doy,\dof)\sq
\Tau_{u\we i,v\we i}(\dox,\doy,\dof)$;

\itla{iii}
if $u,\,v\in 2^n$ and $u(n-1)\not=v(n-1)$ then 
$\Pi_u(\dox)\cup \Pi_v(\doz)\in \dpwo,$ moreover $\in D^2(n),$ 
and $\Pi_u(\dox)\cup \Pi_v(\dox)$ is \dl inconsistent;

\itla{vii}
$\msur\hax{\Tau_{uv}(\dox,\doy,\dof)}=\Pi_u(\dox)$ and 
$\hay{\Tau_{uv}(\dox,\doy,\dof)}=\Pi_v(\doy)$ --- then in 
particular $\hax{\Pi_u(\dox)}=\Pi_u(\dox)$ for all $u$.
\een

\bre
\label R
Since theories in $\dpwt$ contain $\dox\eqh\doy,$ 
it follows from \ref{vii} by Lemma~\ref{dkle} that 
$\hhx{\Pi_u(\dox)}= \hhx{\Pi_v(\dox)}$ for all crucial pairs  
$\ang{u,v}.$ Therefore 
$\hhx{\Pi_u(\dox)}= \hhx{\Pi_v(\dox)}$ for all $u,\,v\in 2^n$ 
and $n\in\om$ as any two tuples $u,\,v\in 2^n$ are connected 
by a (unique) chain of crucial pairs.\qed
\ere

Let us first of all demonstrate that the existence of such a 
system yields a continuous map $F$ in $\kvp$ which witnesses 
\ref2 of Theorem~\ref{h}.

Lemma~\ref{dpa-b} and \ref0 imply that for any $a\in 2^\om$ there 
is a unique real, denoted by $F(a),$ satisfying every formula in 
$\bigcup_{n\in\om}\Pi_{a\res n}(\dox),$ and the map 
$F$ is continuous. Moreover $F$ is $1-1$ by \ref{iii}. 

Suppose that $a,\,b\in 2^\om$ and $a\nEo b,$ so that 
$a(n)\not=b(n)$ for infinitely many $n.$ It follows then from 
\ref{iii} and Theorem~\ref{dpao-b} that $F(a)\nmeks F(b)$. 

Let us check that $F$ satisfies \ref{2} of Theorem~\ref h. 
Suppose that $a,\,b\in 2^\om$ are \dd\meo neighbours, \ie\ 
$a=1^k\we 0\we c$ and $b=0^k\we 1\we c$ for some $k\in\om$ and 
$c\in 2^\om.$ Then $\ang{a\res n,b\res n}$ is a crucial pair for 
all $n>k.$ Therefore, by \ref o and Lemma~\ref{dpat-b}, 
there is a unique triple $\ang{x,y,f}\in\cN^2\ti\rdi^\om$ which 
satisfies every formula in 
$\bigcup_{n\in\om}\Tau_{a\res n,b\res n}(\dox,\doy,\dof),$ and 
now $x=F(a),\msur$ $y=F(b)$ by \ref{vii}. This implies 
$F(a)\mekt F(b)$ by \ref{o}. 

\punk{The construction of theories}

{\it We argue in\/} $\kvp$. 

To define $\Pi_\La$ (where $\La$ is the empty sequence, the 
only member of $2^0$) consider first the theory $\Xi(\dox),$ 
see Subsection~\ref f. As clearly $\Xi(\dox)\in\dP,$ there 
is a theory $\Pi(\dox)\in D(0)$ including $\Xi(\dox).$ 
Let $\Pi_\La(\dox)=\Pi(\dox)$.

Suppose that the construction has been completed up to a level 
$n,$ and expand it to the next level. 

To start with we set $\Pi_{s\we i}(\dox)=\Pi_s(\dox)$ for all 
$s\in 2^n$ and $i=0,1,$ and 
$\Tau_{s\we i,t\we i}(\dox,\doy,\dof)=\Tau_{st}(\dox,\doy,\dof)$ 
whenever $i=0,1$ and $\ang{s,t}$ is a crucial pair in $2^n.$ 
For the ``initial'' pair $\ang{1^n\we 0,0^n\we 1},$ 
let $\Tau_{1^n\we 0,0^n\we 1}$ be the theory  
$\Pi_{0^n}(\dox)\cup\Pi_{0^n}(\doy)\cup(\dox\eqh\doy)\cup 
\ans{\dox\mektf\doy}.$ 
Then, by Lemma~\ref{k2}, $\Tau_{1^n\we 0,0^n\we 1}\in\dpwt$ and 
$\hax{\Tau_{1^n\we 0,0^n\we 1}(\dox,\doy,\dof)}=
\Pi_{1^n\we 0}(\dox),\msur$ 
$\hay{\Tau_{1^n\we 0,0^n\we 1}(\dox,\doy,\dof)}=
\Pi_{0^n\we 1}(\doy)$.

This ends the definition of ``initial values'' at 
the \dd{n\pone}th level. The plan is to gradually enforce the 
theories in order to fulfill the requirements.\vtm

{\bfit Step 1.\/}
Take care of item \ref 0. Consider an arbitrary 
$u_0=s_0\we i\in 2^{n+1}.$ As $D(n)$ is dense there is a theory 
$\Pi'(\dox)\in D(n)$ including $\Pi_{u_0}(\dox).$ We can assume 
that $\Pi'(\dox)=\hax{\Pi'(\dox)}$ for otherwise change 
$\Pi'(\dox)$ by $\hax{\Pi'(\dox)}$. 

The intension is to take $\Pi'(\dox)$ as the 
``new'' $\Pi_{u_0}.$ But this change has to be expanded through 
the net of crucial pairs, in order to preserve 
\ref{vii}. (Fortunately the tree of all crucial pairs in 
$2^{n+1}$ is a chain.) 

Thus put $\Pi'_{u_0}(\dox)=\Pi'(\dox).$ Suppose that 
$\Pi'_u$ has been defined, includes $\Pi_u,$ the older 
version, and satisfies $\hax{\Pi'_u(\dox)}=\Pi'_u(\dox),$ 
for some $u\in 2^{n+1}$ which is connected by a 
crucial pair with a not yet encountered $v\in 2^{n+1}.$ Define 
$\Tau'_{uv}(\dox,\doy,\dof)$ to be 
$\Pi'_u(\dox)\cup \Tau_{uv}(\dox,\doy,\dof)$ and  
$\Pi'_v(\doy)$ to be $\hay{\Tau'_{uv}(\dox,\doy,\dof)}.$ Note 
that $\Pi'_v(\doy)$ includes $\Pi_v(\doy)$ because \ref{vii} 
is assumed for the old theories $\Pi_u,\,\Pi_v,\,\Tau_{uv}.$
Note that \ref{vii} holds for the new theories 
$\Pi'_u,\msur$ $\Pi'_v,\msur$ $\Tau'_{uv}:$ indeed 
$\hax{\Tau'_{uv}(\dox,\doy,\dof)}=\Pi'_u(\dox)$ follows 
from Lemma~\ref{str-o}. 

The construction describes how the change from $\Pi_{u_0}$ to 
$\Pi'_{u_0}$ spreads through the chain of crucial pairs in 
$2^{n+1},$ resulting in a system of new theories, 
$\Pi'_u$ and $\Tau'_{uv},$ which 
satisfy \ref 0 for the particular $u_0\in 2^{n+1}.$ 

We iterate this construction consecutively for all 
$u_0\in 2^{n+1},$ getting finally a system of theories satisfying 
\ref 0 (fully) (and \ref{vii}) which we shall denote by 
$\Pi_u$ and $\Tau_{uv}$ from now on.\vtm

{\bfit Step 2.\/}
Take care of item \ref{iii}. Let us fix a pair of $u_0$ and $v_0$ 
in $2^{n+1},$ such that $u_0(n)=0$ and $v_0(n)=1.$ By the density 
of $D^2(n),$ there is a theory   
$\Pi'_{u_0}(\dox)\cup \Pi'_{v_0}(\doy)\in D^2(n)$ 
which includes  $\Pi_{u_0}(\dox)\cup \Pi_{v_0}(\doy).$ We may 
assume that $\Pi'_{u_0}(\dox)=\hax{\Pi'_{u_0}(\dox)}$ and 
$\Pi'_{v_0}(\doy)=\hay{\Pi'_{v_0}(\doy)}.$ We can also assume, by 
Corollary~\ref{kd}, that 
$\Pi'_{u_0}(\dox)\cup \Pi'_{v_0}(\dox)$ is \dl inconsistent. 

Spread the change from $\Pi_{u_0}$ to $\Pi'_{u_0}$ and from 
$\Pi_{v_0}$ to $\Pi'_{v_0}$ through the chain of crucial pairs in 
$2^{n+1}$ until the two waves of spreading meet each other at the 
pair $\ang{1^n\we 0,0^n\we 1}.$ This leads to a 
system of theories $\Pi'_u$ and $\Tau'_{uv}$ which satisfy 
\ref{iii} for the particular pair $\ang{u_0,v_0}$ and still 
satisfy \ref{vii} with the exception of the  
``meeting'' crucial pair $\ang{1^n\we 0,0^n\we 1}$ (for which 
basically $\Tau'_{1^n\we 0,0^n\we 1}$ is not yet defined for 
this step). 

Take notice that the construction of 
Step 1 has left $\Tau_{1^n\we 0,0^n\we 1}$ in the form 
$\Pi_{1^n\we 0}(\dox)\cup \Pi_{0^n\we 1}(\doy)\cup(\dox\eqh\doy) 
\cup\ans{\dox\mektf\doy}$ 
(where $\Pi_{1^n\we 0}$ and $\Pi_{0^n\we 1}$ are the  
``versions'' at the end of Step 1). We now have new 
\hbox{\dl c}on\-sis\-tent theories, 
$\Pi'_{1^n\we 0}$ and $\Pi'_{0^n\we 1},$ including resp.\ 
$\Pi_{1^n\we 0}$ and $\Pi_{0^n\we 1}$ and satisfying 
$\hhx{\Pi'_{0^n\we 0}(\dox)}=\hhx{\Pi'_{0^n\we 1}(\dox)}.$ 
(See Remark~\ref R; recall that 
$\hhx{\Pi'_{u_0}}=\hhx{\Pi'_{v_0}}$ for the initial pair 
simply because 
$\Pi'_{u_0}(\dox)\cup \Pi'_{v_0}(\doy)\in \dpwo.$) 
Now the theory 
$\Pi'_{1^n\we 0}(\dox)\cup \Pi'_{0^n\we 1}(\doy)
\cup(\dox\eqh\doy)\cup\ans{\dox\mektf\doy}$ 
taken as $\Tau'_{1^n\we 0,0^n\we 1}$ belongs to $\dpwt$ and 
satisfies \ref{vii} for the pair $\ang{1^n\we 0,0^n\we 1}$ 
by Lemma~\ref{k2}. This ends the consideration of the pair 
$\ang{u_0,v_0}$.

Applying this construction consecutively for all pairs of 
$u_0\in P_0$ and $v_0\in P_1$ (including the pair 
$\ang{1^n\we 0,0^n\we 1}$) we finally get a system of theories 
satisfying \ref 0, \ref{iii}, and \ref{vii}, which will be 
denoted still by $\Pi_u$ and $\Tau_{uv}$.\vtm

{\bfit Step 3.\/}
We finally take care of \ref o. Consider a 
particular crucial pair $\ang{u_0,v_0}$ in $2^{n+1}.$ By the 
density of $D_2(n),$ there is a theory   
$\Tau'_{u_0,v_0}(\dox,\doy,\dof)$ in $D_2(n)$ including  
$\Tau_{u_0,v_0}(\dox,\doy,\dof).$ 

Define $\Pi'_{u_0}(\dox)=\hax{\Tau'_{u_0,v_0}(\dox,\doy,\dof)}$ 
and $\Pi'_{v_0}(\doy)=\hay{\Tau'_{u_0,v_0}(\dox,\doy,\dof)}$ and 
spread this change through the chain of crucial pairs in 
$2^{n+1}.$ 
(Note that $\hhx{\Pi'_{u_0}(\dox)}=\hhx{\Pi'_{v_0}(\dox)}$ 
because theories in $\dpwt$ include $\dox\eqh\doy.$ 
This implies $\hhx{\Pi'_u(\dox)}=\hhx{\Pi'_v(\dox)}$ 
for all $u,\,v\in 2^{n+1}$ after the spreading.) 
 
Executing this construction for all crucial pairs 
in $2^{n+1}$ we finally end the construction, in $\kvp,$ of a 
system of theories satisfying \ref 0 through \ref{vii}.\vtm

\qeD{Theorem~\ref h}
\np

\parf{Why embedding $\protect\meo$ is absolute}
\label{abs}

The aim of this section is to prove that \ref2 of Theorem~\ref h 
and \ref{2A} of Theorem~\ref{ha} are absolute statements. By the 
way this fills the gap left in the proof of Theorem~\ref h (see 
the beginning of Section~\ref{E-b}).

\ble
\label{l-abs}
If\/ $p\in\cN$ and\/ $\mek$ is a\/ $\is11(p)$ \pqo\ then\/ 
{\rm\ref{2A}} of Theorem~\ref{ha} is equivalent to a\/ $\is12(p)$ 
statement uniformly in\/ $p.$ Similarly if\/ ${\mekt}\sq{\meks}$ 
are resp.\/ $\is11(p)$ \pqo\ and\/ $\ip11(p)$ \pqo\ then\/ 
{\rm\ref{2}} of Theorem~\ref h is equivalent to a\/ $\is12(p)$ 
statement uniformly in\/ $p$. 
\ele
\bpf
\footnote
{\ The proof involves an idea communicated to the author by 
G.~Hjorth with a reference to Hjorth and Kechris~\cite u, 
Section 3, where the idea is realized in terms of category.
}
We consider only \ref2 of Theorem~\ref h; the other statement is 
pretty similar. The aim does not seem easy: at the first look the 
statement is $\is13.$ To reduce it to $\is12$ we use 
{\it Borel approximations\/} of $\mekt$ and $\meks.$

Recall that 
$\wo=\ans{z\in\cN:z\,\hbox{ codes an ordinal}};$ 
for $z\in\wo$ let $|z|$ be the ordinal coded by $z,$ and 
$\wo_\nu=\ans{z\in\wo:|z|=\nu}$.

Being a $\fs11$ subset of $\cN^2,$ the relation $\mekt$ 
classically has the form  
${\mekt}=\bigcup_{\nu<\omi}{\Rmekt\nu}$ where 
$\ang{{\Rmekt\nu}:\nu<\omi}$ is an increasing sequence of Borel 
subsets of $\mekt.$ 
Moreover there is a $\Pi^1_1$ formula $\pi(z,x,y)$ (containing 
$p$ as a parameter) such that 
$x\Rmekt\nu y\eqv\pi(z,x,y)$ whenever $z\in\wo_\nu.$ (There 
also exists a $\Sg^1_1$ formula with the same property.) 

Similarly 
${\meks}=\bigcap_{\nu<\omi}{\Rmeks\nu},$ where 
$\ang{{\Rmeks\nu}:\nu<\omi}$ is an increasing sequence of Borel 
supersets of $\meks,$ and there is a $\Sg^1_1$ formula 
$\sg(z,x,y)$ such that $x\Rmeks\nu y\eqv\sg(z,x,y)$ whenever 
$z\in\wo_\nu$.

The following statement is clearly $\is12(p)$ (use formulas $\pi$ 
and $\sg$): 
\ben
\def\theenumi{${\mtho({\rm\Roman{enumi}}')}$}
\def\labelenumi{\theenumi}
\addtocounter{enumi}{1}
\itla{2'}
{\it There is a continuous\/ $1-1$ map\/ $F':2^\om\lra\cN$ and a 
countable ordinal\/ $\nu$ such that\/}
\newlength{\ile}
\settowidth{\ile}{(a) \ $a\meo b$}
\newlength{\tle}
\settowidth{\tle}{(b) \ $a\nEo b$}
\newlength{\mle}
\settowidth{\mle}{(a) \ $a\meo b$\ \ {\it implies\/} \ 
$\phantom{\net\;}F'(a)\Rmekt\nu F'(b)\;;$}
\begin{minipage}[t]{\mle}
(a) \ 
$a\meo b$\ \ {\it implies\/} \ 
$\phantom{\net\;}F'(a)\Rmekt\nu F'(b)\;;$\\[1mm]
(b) \
$a\nEo b$\hspace{\ile}\hspace{-\tle}\ \ 
{\it implies\/} \ $\net\; F'(a)\Rmeks\nu F'(b)\;.$ 
\end{minipage}
\hfill$\,$
\een
Thus it remains to prove that \ref2 of Theorem~\ref h is 
equivalent to \ref{2'}. The hard part is to prove that \ref2 
implies \ref{2'}.  
To prove this direction we consider a \dd\kpa collapse generic 
extension $\kvp$ of $\kV,$ the universe of all sets, where $\kpa$ 
is $2^{\aleph_0}$ in $\kV.$ As \ref2 is $\is13(p)$ while \ref{2'} 
is $\is12(p),$ it suffices to prove that \ref2 
implies \ref{2'} {\it in\/ $\kvp.$} 

We can enumerate in $\kvp$ by natural numbers all dense subsets 
of $2\lom$ and $2\lom\ti 2\lom$ (the Cohen forcing and its 
square) which belong to $\kV.$ This allows to define in $\kvp$ 
infinite sequences $\ang{u_n:n\in\om}$ and $\ang{v_n:n\in\om}$ 
such that $u_n,\,v_n\in 2^{l(n)}$ for some $l(n)$ for all $n,$ and 
for any $n$:
\ben
\def\theenumi{$(\fnsymbol{enuF})$}
\def\labelenumi{\theenumi}
\itla,
if $u,\,v\in 2^l$ where $l=n+\sum_{m=0}^{n-1}l(n)$ then the pairs 
$\ang{u\we u_n\,,\,v\we v_n}$ and $\ang{u\we v_n\,,\,v\we u_n}$ 
belong to the \dd nth dense subset of $2\lom\ti 2\lom$.
\addtocounter{enuF}{1}
\een
Define in $\kvp,$ for each $a\in 2^\om,$ 
$G(a)=w_0\we w_1\we w_2\we\dots,$ where $w_n=u_n\we 0$ whenever 
$a(n)=0$ and $w_n=v_n\we 1$ whenever $a(n)=1.$ Then $G$ is 
continuous and $1-1,$ therefore the map $F'(a)=F(G(a))$ is 
continuous and $1-1$ as well. (Here $F$ is a map which witnesses 
\ref2 of 
Theorem~\ref h in $\kvp.$) Prove that $F'$ witnesses \ref{2'}.

Suppose that $a,\,b\in 2^\om$ and $a\meo b.$ Then by definition 
both $a'=G(a)$ and $b'=G(b)$ are {\tt Cohen\/} generic over $\kV$ 
and $a'\meo b',$ $F(a')\mekt F(b')$ (by the choice of $F$), even in 
$\kV[a',b'],$ which implies $F(a')\Rmekt\nu F(b')$ for an 
appropriate $\nu<\omi^{\kV[a',b']}.$ Since the difference between 
$a'$ and $b'$ is finite the latter statement is a 
property of $a',$ hence it is {\tt Cohen\/} forced over $\kV.$ 
It follows that there is a countable in $\kV[a']$ 
(hence in $\kvp$) ordinal $\nu_T$ such that we have, in $\kvp,$ 
$F'(a)\Rmekt{\nu_T} F'(b)$ whenever $a,\,b\in 2^\om$ satisfy 
$a\meo b$.

Suppose that $a,\,b\in 2^\om$ and $a\nEo b.$ Then by definition 
$\ang{G(a),G(b)}$ is ${\tt Cohen}^2$ generic over $\kV,$ in 
particular $a'=G(a)$ and $b'=G(b)$ satisfy $a'\nEo b',$ therefore 
$F(a')\nmeks F(b'),$ which implies $\net\,F(a')\Rmeks\nu F(b')$ 
for an ordinal $\nu<\omi^{\kV[a',b']}.$ As above there is a 
single ordinal $\nu_S<\omi^{\kvp}$ such that we have, in $\kvp,$ 
$\net\,F'(a)\Rmeks{\nu_S} F'(b)$ whenever $a\nEo b$.

It remains to take $\nu={\tt max}\,\ans{\nu_T,\nu_S}$.
\epf

\np

\parf{Borel and analytic order relations}
\label{A}

This section proves theorems \ref{hb} and \ref{ha}, with the 
exception of the additional statement in \ref{1A} of 
Theorem~\ref{ha}. 

\punk{Borel orders: Theorem~\protect\ref{hb}}
\label{bor}

Let $\mek$ be a Borel \pqo\ on $\cN.$ In view of Theorem~\ref h it 
suffices to prove that \ref{b1} and \ref{b2} of Theorem~\ref{hb} 
are incompatible. 
{\it Suppose otherwise\/}. 

The superposition of the maps $F$ and $h$ is then a Borel 
\hop\ map $\phi:\stk{2^\om}{\meo}\lra \stk{2^\al}{\mel}$ 
satisfying the following: $\phi(a)=\phi(b)$ implies that 
$a\Eo b,$ \ie\ $a$ and $b$ are \dd\meo comparable. 

Therefore, as any \dd\Eo class is \dd\meo ordered as $\dZ,\msur$ 
$\om,$ or $\om^\ast,$ the \dd\phi image $X_a=\phi\ima{[a]_{\Eo}}$ 
of the \dd{\Eo}class of any $a\in 2^\om$ is \dd\mel ordered as a 
subset of $\dZ.$ If $X_a=\ans{x_a}$ is a singleton then put 
$\psi(a)=x_a$.  

Assume that $X_a$ contains at least two points. In this case 
we can effectively pick an element in $X_a\,!$ Indeed there is a 
maximal sequence $u\in 2^{<\al}$ such that $u\subset x$ for each 
$x\in X_a.$ Then the set 
$X_a^{\rm left}=\ans{x\in X:u\we 0\subset x}$ 
contains the \dd\mel largest element, which we denote by $\psi(a)$. 

To conclude $\psi$ is a Borel reduction of $\Eo$ to the equality 
on $2^\al,$ \ie\ $a\Eo b$ iff $\psi(a)=\psi(b),$ which is 
impossible.\vtm

\qeD{Theorem~\ref{hb}}

\punk{The general case of analytic relations}
\label{Ag}
\label{rt}

Consider an analytic \pqo\ $\mek$ on $\cN.$ We shall \wlg\ assume 
that $\mek$ is $\is11,$ so that ${\mek}={\mekt}=\pr[T],$ where $T$ 
is 
a recursive tree in $(\om\ti\om\ti\om)\lom.$ We also suppose that
\ben
\def\theenumi{$(\fnsymbol{enuF})$}
\def\labelenumi{\theenumi}
\itla{f1}
$\mekt$ does \underline{not} satisfy \ref{2A} of Theorem~\ref{ha}.
\addtocounter{enuF}{1}
\een
The aim is to prove that then $\mekt$ satisfies \ref{1A} 
(still leaving apart the additional part) of Theorem~\ref{ha}. 

Recall that $\rH_\al$ is defined in Section~\ref{di}. Let 
$\rha=\bigcup_{\al<\omi}\rH_\al\,;$ this includes the set $\rH$ 
also defined in Section~\ref{di}. By definition each $h\in\rha$ 
is a Borel function of certain type, coded in $\kL.$ Let us fix 
a {\it constructible in the codes\/} enumeration 
$\rha=\ans{h_\al:\al<\omi}.$ Define the concatenation
\dm
\gh(x)=h_0(x)\we h_1(x)\we \dots \we h_\al(x) \we \dots 
\hspace{5mm}(\al<\omi)\;.
\dm
Prove that $\gh:\stk\cN\mekt\lra\stk{2\lomi}{\mel}$ is a 
linearization. First of all $\gh$ is a \hop\ map from 
$\stk\cN\mekt$ to $\stk{2^{\omi}}\mel$ because each $h_\al$ is 
\hop\ by definition. Thus it remains to prove that 
$\gh(x)=\gh(y)\imp x\eet y$.

This involves a {\it reflection lemma\/} for analytic \pqo's.

Being a $\fs11$ set, $\mekt$ has the form 
${\mekt}=\bigcap_{\nu<\omi}{\mekta\nu},$ where each $\mekta\nu$ 
is an \dd\omi Borel set coded in\/ $\kL$ 
and ${\mekta\mu}\sq{\mekta\nu}$ provided $\nu<\mu.$ In addition we 
have the following {\bfit boundedness principle\/}: 
if ${\mekt}\sq X,$ where $X\sq\cN^2$ is a $\fp11$ set, then there 
is $\nu<\omi$ such that ${\mekta\nu}\sq X$. 

Take notice that the sets $\mekta\nu$ are not necessarily \pqo's. 
However 

\ble
\label{ref}
Assume that\/ $B\sq\cN^2$ is a Borel set 
and\/ ${\mekt}\sq B.$ Then there is\/ $\nu<\omi$ such that\/ 
${\mekta\nu}\sq B$ and\/ $\mekta\nu$ is a \pqo.~\footnote
{\rm\ Compare with the Burgess reflection theorem for $\fs11$ 
equivalence relations.}
\ele
\bpf
Let us first prove a weaker statement: 
{\it there is\/ $\mu<\omi$ such that\/}
\dm
x\mekta\mu y\,\imp\,\kaz x'\;\kaz y'\;
(x'\mekt x\cj y\mekt y'\imp x'\mekta\mu y')\,.
\dm
Note that by the boundedness there exists an ordinal $\mu_0<\omi$ 
such that ${\mekta{\mu_0}}\sq B.$ Suppose that an ordinal 
$\mu_n\geq \mu_0$ has been defined. Put $C(x,y)$ iff 
$\kaz x'\:\kaz y'\:(x'\mekt x\cj y\mekt y'\imp x'\mekta{\mu_n}y'),$ 
so that $C$ is a $\fp11$ set and ${\mekt}\sq C\sq {\mekta{\mu_n}}.$ 
Using the boundedness principle again we get an ordinal 
$\mu_{n+1}\geq \mu_n$ satisfying 
${\mekt}\sq {\mekta{\mu_{n+1}}}\sq C,$ 
so that by definition 
\dm
x\mekta{\mu_{n+1}}y\,\imp\,\kaz x'\;\kaz y'\;
(x'\mekt x\cj y\mekt y'\imp x'\mekta{\mu_{n}}y')\,.
\dm
It remains to define $\mu=\sup_n \mu_n$. 

Starting the proof of the lemma, we choose $\nu_0$ so that 
${\mekta{\nu_0}}\sq B$ and 
\dm
x\mekta{\nu_0}y\,\imp\,\kaz x'\;\kaz y'\;
(x'\mekt x\cj y\mekt y'\imp x'\mekta{\nu_0}y')\,.
\eqno{(\ast)}
\dm
Suppose that an ordinal $\nu_n\geq \nu_0$ satisfying $(\ast)$ has 
been defined. Put $C(x,y)$ iff 
$x\mekta{\nu_n}y \cj \kaz z\:(y\mekta{\nu_n}z\imp x\mekta{\nu_n}z),$ 
so that $C$ is a $\fp11$ set and ${\mekt}\sq C\sq {\mekta{\nu_n}}.$ 
As above there is an ordinal $\nu_{n+1}\geq\nu_n$ satisfying 
${\mekt}\sq {\mekta{\nu_{n+1}}}\sq C,$ so that by definition 
\dm
x\mekta{\nu_{n+1}}y\,\imp\,
\kaz z\:(y\mekta{\nu_n}z\imp x\mekta{\nu_n}z)\,.
\dm
Now take $\nu=\sup_n\nu_n$.
\epf

Now assume $x\not\eet y$ and prove that $\gh(x)\not=\gh(y).$ 
By the lemma there is an ordinal $\nu<\omi$ such that 
$x\not\eet^\nu y$ and $\mekta\nu$ is a \pqo. 
It is a classical fact that the {\it complement\/} of $\mekta\nu$ 
(which is the \dd\nu th {\it constituent\/} 
of the co-analytic set $\dop({\mekt})$) 
can be presented in the form of the set of all pairs $\ang{x,y}$ 
such that $C(x,y)$ is well-ordered and has the order type 
$\<\nu,$ where $C:\cN^2\lra\dQ$ is a continuous function coded in 
$\kL[T]=\kL.$ Therefore there exists a tree $S=S_\nu\in\kL,\msur$ 
$S\sq(\om\ti\om\ti\nu)\lom,$ such that $\mekta\nu=\dop\pr[S].$ 

Apply Theorem~\ref h for the relations 
${\mekt}=\pr[T]\sq{\mekt^\nu}=\dop\pr[S].$ 
We observe that \ref2 of Theorem~\ref h fails by the assumption 
\ref{f1}. Therefore \ref1 of Theorem~\ref h holds, 
so that there exists an ordinal $\al<\omi$ such that the map 
$h_\al$ satisfies $h_\al(x)=h_\al(y)\imp x\eet^\nu y$. 

It follows that $h_\al(x)\not=h_\al(y)$ by the choice of $\nu,$ 
hence $\gh(x)\not=\gh(y),$ as required.
Thus $\gh$ witnesses \ref{1A} (the general case) of 
Theorem~\ref{ha}.

\np

\parf{Special cases: Borel classes and generic models}
\label{S}

This section is devoted to the ``additional'' part in \ref{1A} 
of Theorem~\ref{ha}. 
{\it We continue to argue in the assumptions and notation\/ 
of Subsection~\ref{Ag}.\/} 
In particular we still assume \ref{f1} of Subsection~\ref{Ag}, 
that is the given $\is11$ \pqo\ $\mekt=\pr[T]$ does {\it not\/} 
satisfy \ref{2A} of Theorem~\ref{ha}.


We need here to specify definability properties related to 
the function $\gh.$ An ordinary argument allows to 
define the sequence ${\ang{h_\al:\al<\omi}}$ to be 
$\id{\kL_{\omi}}1$ hence $\id\hc1.$ Then the map 
$\ang{x,\ga}\mapsto\gh(x)\res\ga$ is $\id\hc1,$ too. 
Moreover there is a $\Sg_1$ formula $\Phi(\cdot,\cdot,\cdot)$ 
such that 
\ben
\def\theenumi{$(\fnsymbol{enuF})$}
\def\labelenumi{\theenumi}
\itla{Phi}
If $M$ is a transitive model of $\ZFC^-$ (minus Power Sets), 
$x\in\cN\cap M,$ and $\la\in M,\msur$ $\la<\omi,$ then 
$u=\gh(x)\res\la\in M$ and $u$ is the only member of $M$ such 
that $\Phi(x,\la,u)$ holds in $M$.
\addtocounter{enuF}{1}
\een

\punk{Order relations with Borel classes}
\label{B}

Consider the case \ref{1Abor} in the ``additional'' part of 
Theorem~\ref{ha}. Suppose that 
{\it every\/ \dd\eet class\/ $[x]_{\eet}=\ans{y:y\eet x}$ is 
Borel\/.} The aim is to find an antichain $A\sq 2\lomi$ and a 
$\id\hc1$ linearization 
$\gh':\stk\cN\mekt\lra \stk A\mel$.\vtm

{\bfit First attempt.\/} 
Let $x\in\cN.$ As $[x]_{\eet}$ is Borel, there exists, by 
Lemma~\ref{ref}, an ordinal $\la<\omi$ such that 
$x\eet y\eqv x\eet^\la y$ and $\mekt^\la$ is a \pqo. Therefore 
(see the end of Subsection~\ref{Ag}) there is an ordinal 
$\al<\omi$ satisfying $h_\al(x)=h_\al(y)\imp x\eet^\nu y,$ hence 
$h_\al(x)=h_\al(y)\eqv x\eet y$ (for the chosen $x$ and any $y$). 
It follows that $x\eet y\eqv \gh(x)\res\la = \gh(y)\res\la,$ for 
an ordinal $\la<\omi.$ Let $\la_0(x)$ be the least such an ordinal 
and $\gh_0(x)=\gh(x)\res\la_0(x).$ Now 
$A=\ans{\gh_0(x):x\in\cN}$ is an antichain in $2^{<\omi}$ while 
$\gh_0:\stk\cN\mekt\lra\stk A\mel$ is a linearization. 

However the definition of $\la_0(x)$ seems not to provide that 
$\gh_0$ is $\id\hc1$.\vtm

{\bfit Second attempt.\/} Let us modify the construction to meet 
the requirement that $\gh_0$ is $\id\hc1.$ 
The idea is to replace the quantifier $\kaz y$ in the definition 
of $\la_0(x)$ by: 
{\it for each\/ $y$ in every set generic extension of 
a suitable model containing\/ $x,$} which essentially means: 
{\it for comeager-many reals~$y$}.

Let $x\in\cN.$ Define $\la_1(x)$ to be the least ordinal $\la<\omi$ 
such that 
\ben
\def\theenumi{$(\fnsymbol{enuF})$}
\def\labelenumi{\theenumi}
\itla{dag}
there exists a transitive model $M$ of $\ZFC^-$ containing $\la,$ 
and a real $x'\in M$ satisfying $x\eet x'$ and, for any set 
generic extension $M'$ of $M$ and each real $z\in M',$ we have 
$\gh(x')\res\la=\gh(z)\res\la\imp x'\eet z$.
\addtocounter{enuF}{1}
\een
(Clearly $\la_1(x)\leq\la_0(x)$.) 

Define $\gh_1(x)=\gh(x)\res\la_1(x).$ Thus $\gh_1:\cN\lra 2\lomi$ 
is a $\id\hc1$ function. 

\ble
\label{oK}
We have$:$ $x\eet y$ iff\/ $\gh_1(x)=\gh_1(y)$.
\ele
\bpf
Suppose that $\gh_1(x)=\gh_1(y)=u\in 2\lomi$ 
(the nontrivial direction), in particular 
$\la_1(x)=\la_1(y)=\la.$ Let $x'\in M_x$ and $y'\in M_y$ witness 
that $\la_1(x)=\la_1(y)=\la,$ in the sense of \ref{dag}. Prove 
that $x'\eet y'.$ 

Let $f\in\la^\om$ be a \dd\la collapse function generic over both 
$M_x$ and $M_y.$ Let $\vt_x,\,\vt_y$ be the least ordinals not in 
resp.\ $M_x,\,M_y.$ Suppose that $\vt_x\leq\vt_y.$ Then both 
$M_x[f]$ and $\kL_{\vt_x}[f]\sq M_x[f]$ model $\ZFC^-,$ hence by 
some version of Shoenfield there is a real $z\in\kL_{\vt_x}[f]$ 
satisfying $\gh(z)\res\la=u,$ therefore $x'\eet z$ by the assumed 
\ref{dag}. However $z$ can be 
also considered in $\kL_{\vt_y}[f]$ which yields $y'\eet z,$ 
hence $x'\eet y'$ as required.
\epf

Thus we have defined a $\id\hc1$ map $\gh_1:\cN\lra 2\lomi$ 
such that, for any $x,$ $\gh_1(x)=\gh(x)\res\la_1(x)$ for some 
$\la_1(x)<\omi,$ and ${x\eet y}\eqv{\gh_1(x)=\gh_1(y)}.$ However 
the range $A_1=\ans{\gh_1(x):x\in\cN}$ may be not an antichain in 
$2\lomi.$ 

To fix this last problem, we define a new $\id\hc1$ map 
$\gh':\cN\lra5\lomi.$ First we change all values $\gh_1(x)(\ga)=1$ 
to $\gh'(x)(\ga)=4$ and add that $\gh'(x)(\la_1(x))=2.$ If 
$\nu<\la_1(x)$ is an ordinal of the form $\nu=\la_1(y)$ for 
some $y$ then we change $\gh'(x)(\nu)$ once again: if $x\mekt y$ 
then put $\gh'(x)(\nu)=1$ while if $y\mekt x$ then put 
$\gh'(x)(\nu)=3.$ 
A simple verification shows that $\gh'$ satisfies \ref{1A} of 
Theorem~\ref{ha} (except for the fact that $\gh'$ takes values in 
$5\lomi$ rather than $2\lomi$ but this can be easily fixed).

\punk{Order relations in a generic universe}
\label{G}

This subsection is devoted to subitem~\ref{1Agen} in 
Theorem~\ref{ha}. In fact we shall assume the following: 
{\it there exists a real\/ $z_0$ such that each real\/ $x$ in the 
universe\/ $\kV$ belongs to a set~\footnote
{\ It is not clear to what extent {\it class\/} forcing universes 
can accomodate the reasoning below, in particular the proofs of 
lemmas \ref{kly} and \ref{eta}. 
} 
generic extension of\/ $\kL[z_0].$}~\footnote
{\ The extensions can be different for different reals $x.$ 
Moreover the extensions can be Boolean valued extensions of 
$\kL[z]$ rather than factual classes in the universe.}

It can be assumed that in fact $z_0=0,$ so we simply drop $z_0$.

\ble
\label{kly}
Let\/ $x\in\cN.$ There is an ordinal\/ $\la<\omi$ such that\/ 
\ben
\def\theenumi{{\protect\rm(\roman{enumi})}}
\def\labelenumi{\theenumi}
\itla i\msur
$\kL_\la[\gh(x)\res\la]$ models\/ $\ZFC^-$\/ 
{\rm(minus Power Sets)}$;$

\itla{ii}\msur
$x$ belongs to a\/ 
set generic extension of\/ $\kL_\la[\gh(x)\res\la]\;;$

\itla{iiI}
if\/ $y,\,z$ are reals in a set generic extension\/ $M$ of\/ 
$\kL_\la[\gh(x)\res\la],$ and\/ 
$\gh(y)\res\la=\gh(z)\res\la,$ then\/ $y\eet z$.
\een
\ele
\bpf
There exists a cardinal $\kpa$ such that $x$ belongs to a set 
generic extension of $\kL_\kpp$ where $\kpp$ is taken in the 
sense of $\kL.$ Note that $\kL_\kpp$ models $\ZFC^-.$ There is 
an ordinal $\la<\omi$ such that $\stk{\kL_\la}{x,{\in}}$ is 
\dd\in isomorphic to a countable elementary submodel of 
$\stk{\kL_\kpp}{x,{\in}},$ hence $\kL_\la$ models $\ZFC^-$ 
and $x$ belongs to a 
set generic extension of $\kL_\la.$ Note that $u=\gh(x)\res\la$ 
is a class in $\kL_\la[x]$ by \ref{Phi}. Now it is a known fact 
(see Lemma 4.4 in Solovay~\cite{sol} or Lemma 5 in 
Kanovei~\cite{k-gen} as particular cases) 
that $x$ belongs to a set generic extension of $\kL_\la[u]$.

Prove \ref{iiI}. Suppose that reals $y,\,z$ belong to a set 
generic extension of $\kL_\la[u]$ and 
$\gh(y)\res\la=\gh(y)\res\la\,;$ prove $y\eet z.$ First of all, 
a standard forcing argument shows that, as $u=\gh(x)\res\la$ is 
definable in 
$\kL_\la[x]$ while $x$ belongs to a set generic extension of 
$\kL_\la$ by the above, we can \wlg\ assume that $y,\,z$ belong 
to a set generic extension of $\kL_\la$ itself.

Let, for a transitive model $M$ of $\ZFC^-,$ $\gh^M$ denote the 
map $\gh$ defined in $M$ 
(via \eg\ the formula $\Phi$ of \ref{Phi}). 
For instance $\gh^\kV$ is simply $\gh.$ In any case $\gh^M$ 
maps reals in $M$ into $2^{{\omi^M}}.$ In particular if each 
ordinal $\al\in M$ is countable in $M$ and $\la=\Ord\cap M$ then 
$\gh^M(x)\in 2^\la$ for any real $x\in M$. 

Note that \ref{2A} of Theorem~\ref{ha} fails in any generic 
extension $\kL[G]$ of $\kL$ because it is essentially a $\is12$ 
sentence (by Lemma~\ref{l-abs}) false in $\kV$ by \ref{f1}. 
Therefore, by the already proved, in Subsection~\ref{Ag}, part 
of Theorem~\ref{ha}, every set forcing notion in $\kL$ forces 
that, in the extension $\kL[G],$ for any two reals $y,\,z,$ 
if $\gh^{\kL[G]}(y)=\gh^{\kL[G]}(z)$ then $y\eet z$. 

A simple forcing argument transfers this result to $\kL_\kpp$ as 
the initial model, and then, by the choice of $\la,$ to $\kL_\la,$ 
so that for any two reals $y,\,z$ in a set generic extension 
$\kL_\la[G]$ of $\kL_\la,$ if 
$\gh^{\kL_\la[G]}(y)=\gh^{\kL_\la[G]}(z)$ then $y\eet z.$ It 
remains to note that 
$\gh^{\kL_\la[G]}(y)=\gh(y)\res\omi^{\kL_\la[G]}$ for any real 
$y\in\kL_\la[G]$. 
\epf

Let, for any $x\in\cN,$ $\la_x$ be the least ordinal $\la<\omi$ 
satisfying the requirements of the the lemma. We put 
$\gh_1(x)=\gh(x)\res\la_x.$ Apparently $\gh_1$ is a $\id\hc1$ map. 

\ble
\label{eta}
We have$:$ $\gh_1(x)=\gh_1(y)$ iff\/ $x\eet y$.
\ele
\bpf
Suppose that $\gh_1(x)=\gh_1(y)=u\in 2^\la$ (so that 
$\la=\la_x=\la_y$) and prove $x\eet y$ 
(the nontrivial direction). Let $x'$ and $y'$ witness that 
$\gh_1(x)=\gh_1(y)=u,$ so that they belong to resp.\ 
$\kL_\la[u,G_x]$ and $\kL_\la[u,G_y],$ which are resp.\ 
\dd{P_x}generic and \dd{P_y}generic extensions of 
$\kL_\la[u]\,;$ $P_x$ and $P_y$ being set forcing notions in 
$\kL_\la[u].$ In addition, $x\eet x'$ and $y\eet y'.$ 

In particular, as $x\eet x',$ we have $\gh(x)=\gh(x'),$ so that 
$u=\gh(x')\res\la.$ Recall that the map $\gh(\cdot)\res\la$ 
results in some effective way from the $\id{\kL_\la}1$ sequence 
$\ang{h_\al:\al<\la}.$ Therefore the fact that $\gh(x')\res\la=u$ 
is forced in $\kL_\la[u].$ Thus we can assume that $P_x$ forces in 
$\kL_\la[u]$ that $\gh(x')\res\la=u$. 

Consider a set $G\sq P_x$ which is \dd{P_x}generic over 
both $\kL_\la[u,G_x]$ and $\kL_\la[u,G_y].$ Let 
$z\in\kL_\la[u,G]$ be produced by $G$ as $x'$ from $G_x,$ 
so that $\gh(z)\res\la=u.$ 
Thus two reals, $x'$ and $z,$ is the model $\kL_\la[u,G_x,G],$ 
satisfy $\gh(z)\res\la=\gh(x')\res\la.$ 
It follows that $x'\eet z$ by the choice of $\la$.

We similarly prove $y'\eet z,$ as required.
\epf

It remains to get $\gh'$ from $\gh_1$ as in Subsection~\ref{B}.

\qeD{Theorem~\ref{ha}\/}

\np

\let\section=\subsection


\small

\begin{thebibliography}{99}
\label{refe}

\bibitem{f}
S. D. Friedman and B. Velickovic. 
Nonstandard models and analytic equivalence relations. 
{\it Proc.\ Amer.\ Math.\ Soc.\/}, 
to appear.

\bibitem{hkl} 
L.\ A.\ Harrington, A.\ S.\ Kechris, A.\ Louveau. 
A Glimm -- Effros dichotomy for Borel equivalence relations.   
{\it J. Amer. Math. Soc.\/} 
1990, 3, pp.\ 903 -- 928.

\bibitem{hms}
L.\ A.\ Harrington, D.\ Marker, S.\ Shelah,
Borel orderings. 
{\it Trans. Amer. Math. Soc.\/} 
1988, 310, pp.\ 293 -- 302.

\bibitem{hs} L.\ A.\ Harrington, S.\ Shelah, 
Counting equivalence classes for 
co-$\mtho\kpa$-Souslin equivalence relations. \ 
D.\ van Dalen \ea\ eds, 
{\it Proceedings of the Conference in Prague\/} (LC'80), 
North Holland, 1982, pp.\ 147 -- 152.

\bibitem{h94}
G.\ Hjorth. 
{\it A remark on\/ $\fp11$ equivalence relations\/}. 
Note.

\bibitem{h-thin}
G.\ Hjorth. 
Thin equivalence relations and effective decompositions. 
{\em J. Symbolic Logic\/} 
1993, 58, pp.\ 1153 -- 1164. 

\bibitem{u}
G.\ Hjorth, A.\ S.\ Kechris. 
Analytic equivalence relations and Ulm--type classification. 
{\it J. Symbolic Logic\/} 
1995, 60, pp.\ 1273 -- 1299.

\bibitem{k-sm} 
V.\ Kanovei. 
An Ulm--type classification theorem for equivalence relations in 
Solovay model. 
{\it J.\ Symbolic Logic\/} 
1997, to appear.

\bibitem{k-g} 
V. Kanovei. 
Two dichotomy theorems on colourability of non-analytic graphs. 
{\it Fund. Math.\/}, 
to appear. 

\bibitem{k-gen}
V. Kanovei.
Ulm classification of analytic equivalence relations in generic 
universes. 
{\it Math.\ Logic Quarterly\/}, 
to appear. 

\bibitem{k-b}
V. Kanovei.
When a partial Borel order is linearizable.
{\it Fund. Math.\/},
submitted.

\bibitem{l}
A.\ Louveau. 
Two results on Borel orders. 
{\it J. Symbolic Logic\/} 
1989, 54, pp.\ 865 -- 874.

\bibitem{s}
S. Shelah.
On co-$\mathsurround=0mm\kappa$-Souslin relations. 
{\it Israel J.\ Math\/}
1984, 47, pp.\ 139 -- 153. 

\bibitem{sol} 
R.\ M.\ Solovay. 
A model of set theory in which every set of reals 
is Lebesgue measurable. 
{\it Ann. Math.\/} 
1970, 92, pp.\ 1 -- 56.

\end{thebibliography}
\end{document}